# RELAXED LYAPUNOV CRITERIA FOR ROBUST GLOBAL STABILIZATION OF NONLINEAR SYSTEMS

**Iasson Karafyllis*, Costas Kravaris** and Nicolas Kalogerakis***

***Department of Environmental Engineering, Technical University of Crete, 73100, Chania, Greece**

****Department of Chemical Engineering, University of Patras 1 Karatheodory Str., 26500 Patras, Greece**

**Abstract**
The notion of relaxed Robust Control Lyapunov Function (RCLF) is introduced and is exploited for the design of robust feedback stabilizers for nonlinear systems. Particularly, it is shown for systems with input constraints that "relaxed" RCLFs can be easily obtained, while RCLFs are not available. Moreover, it is shown that the use of "relaxed" RCLFs usually results to different feedback designs from the ones obtained by the use of the standard RCLF methodology. Using the "relaxed" RCLFs feedback design methodology, a simple controller that guarantees robust global stabilization of a perturbed chemostat model is provided.

**Keywords:** uniform robust global feedback stabilization, chemostat models.

## 1. Introduction

Consider a finite-dimensional control system:

$$\begin{aligned} &\dot{x} = f(d,x) + g(d,x)u \\ &x \in \Re^n \,, d \in D, u \in U \end{aligned} \tag{1.1}$$

where $D \subset \Re^l$ is a compact set, $U \subseteq \Re^m$ a non-empty convex set with $0 \in U$, $f : D \times \Re^n \to \Re^n$, $g : D \times \Re^n \to \Re^{n \times m}$ are continuous mappings with $f(d,0) = 0$ for all $d \in D$. The problem of existence and design of a continuous feedback law $k : \Re^n \to U$ with $g(d,0)k(0) = 0$ for all $d \in D$, which achieves robust global stabilization of $0 \in \Re^n$ for (1.1), i.e., $0 \in \Re^n$ is uniformly robustly globally asymptotically stable for the closed-loop system $\dot{x} = f(d,x) + g(d,x)k(x)$, is closely related to the existence of a Robust Control Lyapunov Function (RCLF) for (1.1), i.e., the existence of a continuously differentiable, positive definite and radially unbounded function $V : \Re^n \to \Re^+$ with

$$\inf_{u \in U} \sup_{d \in D} \nabla V(x)\big(f(d,x) + g(d,x)u\big) < 0 \text{, for all } x \neq 0 \text{, } x \in \Re^n \tag{1.2}$$

The reader should consult [2,4,9,12,14,15,25,26] and references therein, where the methodology of Lyapunov feedback (re)design is explained in detail. However, in many cases it is very difficult to obtain a CLF for a given control system. The goal of the present work is to show that continuously differentiable, positive definite and radially unbounded functions $V : \Re^n \to \Re^+$ with

$$\inf_{u \in U} \sup_{d \in D} \nabla V(x)\big(f(d,x) + g(d,x)u\big) < 0 \text{, for all } x \neq 0 \text{, } x \in \Omega \tag{1.3}$$

where $\Omega \subseteq \Re^n$ does not necessarily coincide with the whole state space $\Re^n$, can be used in order to design a globally stabilizing feedback. Particularly, under appropriate hypotheses, we show that a continuous feedback law $k : \Re^n \to U$ with $g(d,0)k(0) = 0$ for all $d \in D$, which guarantees

$$\sup_{d\in D} \nabla V(x)\big(f(d,x)+g(d,x)k(x)\big)<0 \text{, for all } x\neq 0,\ x\in\Omega \tag{1.4}$$

and for which $\Omega \subseteq \Re^n$ is an absorbing set for the closed-loop system (1.1) with $u=k(x)$, i.e., every solution of the closed-loop system (1.1) with $u=k(x)$ enters $\Omega \subseteq \Re^n$ in finite time, achieves robust global stabilization of $0\in\Re^n$ for (1.1) (see Theorem 2.2 and Theorem 2.6 below). The reader should compare condition (1.3) with condition (1.2): it is much easier task to find continuously differentiable, positive definite and radially unbounded functions $V:\Re^n \to \Re^+$ satisfying (1.3) instead of (1.2). For this reason we will call a continuously differentiable, positive definite and radially unbounded function $V:\Re^n \to \Re^+$ satisfying (1.3) a "relaxed" RCLF.

It should be emphasized that the idea explained above is intuitive and has been used in the literature in one form or another: for example, it has been used (without explicit statement of the idea) in [12] and in [32] for special classes of control systems. In the present work, a general theoretical formulation and results are developed, so as to provide systematic guidelines for the construction of feedback based on a "relaxed" RCLF. The development of the notion of the "relaxed" RCLF leads to two important applications:

i) even if a RCLF is known then the use of the "relaxed" RCLF feedback design methodology usually results to different feedback designs from the ones obtained by the use of the standard RCLF design methodology; particularly, there is no need to make the derivative of RCLF negative everywhere,

ii) in many cases "relaxed" RCLFs can be found, while RCLFs are not available. Consequently, the class of systems where Lyapunov-based feedback design principles can be applied is enlarged.

In order to illustrate the above applications of the notion of "relaxed" RCLF we consider the following applications:

i) The obtained results are applied to the problem of robust feedback stabilization of the chemostat (Section 3), which has recently attracted attention (see [1,6,10,11,13,17,21,22,23] as well as [8,24,33,34] for studies of the dynamics of chemostat models). In this work we consider the robust global feedback stabilization problem for the more general uncertain chemostat model

$$\begin{aligned}
&\dot{X}=\big(\mu(S)+\Delta(S,t)-D-b\big)X\\
&\dot{S}=D(S_i-S)-K\mu(S)X+m\,X\\
&X\in(0,+\infty)\,,S\in(0,S_i)\,,D\geq 0
\end{aligned} \tag{1.5}$$

where $\Delta(S,t)$ represents a vanishing perturbation (uncertainty). The chemostat model (1.5) and the form of the perturbation $\Delta(S,t)$ are explained in Section 3. As far as we know, this is the first time that the robust global feedback stabilization problem for the chemostat model (1.5) is studied and the proposed controllers in the literature cannot guarantee Robust Global Asymptotic Stability for the resulting closed-loop system (for detailed explanations see Section 3 below). Under mild hypotheses for the equilibrium point $(X_s,S_s)$ of system (1.5) a RCLF for (1.5) is given in Proposition 3.1. However, using the standard RCLF approach we obtain very complicated stabilizing feedback laws. Different families of robust global stabilizers are obtained by exploiting the idea of "relaxed" RCLFs. Particularly, we show that for every locally Lipschitz non-increasing function $\psi:(0,S_i)\to\Re^+$ with $\psi(S)=0$ for all $S\geq S_s$ and $\psi(S)>0$ for all $S<S_s$ and for every locally Lipschitz function $L:(0,+\infty)\times(0,S_i)\to(0,+\infty)$ with $\inf\{L(X,S):(X,S)\in(0,+\infty)\times(0,S_i)\}>0$, the locally Lipschitz feedback law:

$$D=\frac{S_s}{S_i-S_s}\max\big(0\,,K\mu(S)-m\big)\frac{X}{S}+L(X,S)\psi(S) \tag{1.6}$$

guarantees robust global asymptotic stabilization of the equilibrium point $(X_s,S_s)$ of system (1.5).

ii) The obtained results are applied to the problem of feedback stabilization of affine in the control nonlinear systems of the form (1.1) with input constraints. The feedback stabilization problem for nonlinear systems with input constraints has attracted attention (see [18,19,20,28,29,30,32]). Using the idea of "relaxed" RCLFs we are able to reproduce (and slightly generalize) the main results in [32] concerning triangular systems with input constraints (Theorem 4.4 below) as well as obtain simple sufficient conditions for the existence of stabilizing feedback with a simple saturation (Example 4.1 below). It is shown for systems with input constraints that "relaxed" RCLFs can be easily obtained, while RCLFs are not available.

**Notations** Throughout this paper we adopt the following notations:

* For a vector $x \in \Re^n$ we denote by $|x|$ its usual Euclidean norm and by $x'$ its transpose.
* We say that an increasing continuous function $\gamma : \Re^+ \to \Re^+$ is of class $K$ if $\gamma(0) = 0$. By $KL$ we denote the set of all continuous functions $\sigma = \sigma(s,t) : \Re^+ \times \Re^+ \to \Re^+$ with the properties: (i) for each $t \ge 0$ the mapping $\sigma(\cdot, t)$ is of class $K$ ; (ii) for each $s \ge 0$, the mapping $\sigma(s, \cdot)$ is non-increasing with $\lim_{t \to +\infty} \sigma(s,t) = 0$.
* Let $D \subseteq \Re^l$ be a non-empty set. By $M_D$ we denote the class of all Lebesgue measurable and locally essentially bounded mappings $d : \Re^+ \to D$.
* By $C^j(A)$ ( $C^j(A\,;\Omega)$ ), where $j \ge 0$ is a non-negative integer, $A \subseteq \Re^n$, we denote the class of functions (taking values in $\Omega \subseteq \Re^m$) that have continuous derivatives of order $j$ on $A$.
* For every scalar continuously differentiable function $V : \Re^n \to \Re$, $\nabla V(x)$ denotes the gradient of $V$ at $x \in \Re^n$, i.e., $\nabla V(x) = \left( \frac{\partial V}{\partial x_1}(x), \ldots, \frac{\partial V}{\partial x_n}(x) \right)$. We say that a function $V : \Re^n \to \Re^+$ is positive definite if $V(x) > 0$ for all $x \ne 0$ and $V(0) = 0$. We say that a continuous function $V : \Re^n \to \Re^+$ is radially unbounded if the following property holds: "for every $M > 0$ the set $\{ x \in \Re^n : V(x) \le M \}$ is compact".
* The saturation function $\Re \ni x \to sat(x)$ is defined by $sat(x) := x$ for $x \in [-1,1]$ and $sat(x) := \frac{x}{|x|}$ for $x \notin [-1,1]$.

## 2. Main Results

In this section the main results of the present work are presented. We start by recalling the notion of Uniform Robust Global Asymptotic Stability. Consider the following dynamical system:

$$\begin{aligned} &\dot{x} = F(d,x) \\ &x \in \Re^n\,, d \in D \end{aligned} \tag{2.1}$$

We assume throughout this section that system (2.1) satisfies the following hypotheses:

**(H1)** $D \subset \Re^l$ is compact.

**(H2)** The mapping $D \times \Re^n \ni (d,x) \to F(d,x) \in \Re^n$ is continuous.

**(H3)** There exists a symmetric positive definite matrix $P \in \Re^{n \times n}$ such that for every compact set $S \subset \Re^n$ it holds that $\sup \left\{ \frac{(x-y)' P (F(d,x) - F(d,y))}{|x-y|^2} : d \in D, x, y \in S, x \ne y \right\} < +\infty$.

Hypothesis (H2) is a standard continuity hypothesis and hypothesis (H3) is often used in the literature instead of the usual local Lipschitz hypothesis for various purposes and is called a "one-sided Lipschitz condition" (see, for example [28], page 416 and [7], page 106). Notice that the "one-sided Lipschitz condition" is weaker than the hypothesis of local Lipschitz continuity of the vector field $F(d,x)$ with respect to $x \in \Re^n$. It is clear that hypothesis (H3) guarantees that for every $(x_0, d) \in \Re^n \times M_D$, there exists a unique solution $x(t)$ of (2.1) with initial condition $x(0) = x_0$ corresponding to input $d \in M_D$. We denote by $x(t; x_0, d)$ the unique solution of (2.1) with initial condition $x(0) = x_0 \in \Re^n$ corresponding to input $d \in M_D$.

**Definition 2.1:** *We say that $0 \in \Re^n$ is uniformly robustly globally asymptotically stable (URGAS) for (2.1) under hypotheses (H1-3) with $F(d,0) = 0$ for all $d \in D$ if the following properties hold:*

- *for every $s > 0$, it holds that*

$$\sup\left\{ |x(t;x_0,d)| \,;\, t \geq 0 \,,\; |x_0| \leq s \;\;, d \in M_D \right\} < +\infty$$

**(Uniform Robust Lagrange Stability)**

- *for every* $\varepsilon > 0$ *there exists a* $\delta := \delta(\varepsilon) > 0$ *such that:*

$$\sup\left\{ |x(t;x_0,d)| \,;\, t \geq 0 \,,\; |x_0| \leq \delta \;\;, d \in M_D \right\} \leq \varepsilon$$

**(Uniform Robust Lyapunov Stability)**

- *for every* $\varepsilon > 0$ *and* $s \geq 0$*, there exists a* $\tau := \tau(\varepsilon, s) \geq 0$*, such that:*

$$\sup\left\{ |x(t;x_0,d)| \,;\, t \geq \tau \,,\; |x_0| \leq s \;\;, d \in M_D \right\} \leq \varepsilon$$

**(Uniform Attractivity for bounded sets of initial states)**

It should be noted that the notion of uniform robust global asymptotic stability coincides with the notion of uniform robust global asymptotic stability presented in [16].

Next we present relaxed Lyapunov-like sufficient conditions for URGAS. The Lyapunov-like conditions of the following theorem are "relaxed" in the sense that the Lyapunov differential inequality is not required to hold for every non-zero state, but only for states that belong to an appropriate set of the state space. On the other hand an additional reachability condition must hold. Its proof is provided in the Appendix.

**Theorem 2.2:** *Consider system (2.1) under hypotheses (H1-3) with* $F(d,0) = 0$ *for all* $d \in D$ *and suppose that there exists a set* $\Omega \subseteq \Re^n$ *with* $0 \in \Omega$*, functions* $V \in C^1(\Omega;\Re^+)$ *being positive definite and radially unbounded,* $T \in C^0(\Re^n;\Re^+)$*,* $G \in C^0(\Re^n;\Re^+)$*, which satisfy the following properties:*

**(P1)** *For every* $(d, x_0) \in M_D \times \Re^n$*, there exists* $\hat{t}(x_0,d) \in [0, T(x_0)]$ *such that the unique solution* $x(t;x_0,d)$ *of (2.1) satisfies* $x(t;x_0,d) \in \Omega$ *for all* $t \in [\hat{t}(x_0,d), t_{\max})$ *and* $|x(t;x_0,d)| \leq G(x_0)$ *for all* $t \in [0, \hat{t}(x_0,d)]$*, where* $t_{\max} = t_{\max}(x_0,d)$ *is the maximal existence time of the solution,*

**(P2)** $\sup_{d \in D}\left(\nabla V(x) F(d,x)\right) < 0$ *for all* $x \in \Omega$*,* $x \neq 0$.

*Then* $0 \in \Re^n$ *is URGAS for (2.1).*

**Remark 2.3:** For disturbance-free systems, hypothesis (P1) of Theorem 2.2 guarantees that the set $\Omega \subseteq \Re^n$ is an absorbing set (see [31]). Notice that the set $\Omega \subseteq \Re^n$ is not required to be positively invariant. In [12] the name "capturing region" was given for the more general case of set-valued maps instead of sets $\Omega \subseteq \Re^n$ having properties (P1) and (P2) for general time-varying systems.

The proof of Theorem 2.2 utilizes the following lemma, which shows that Uniform Robust Lagrange Stability and Uniform Attractivity for bounded sets of initial states are sufficient conditions for URGAS. Its proof is provided in the Appendix.

**Lemma 2.4:** *Consider (2.1) under hypotheses (H1-3) with* $F(d,0) = 0$ *for all* $d \in D$ *and suppose that there exists a continuous function* $R : \Re^n \to \Re^+$ *such that for every* $(x_0, d) \in \Re^n \times M_D$ *the solution* $x(t;x_0,d)$ *of (2.1) satisfies for all* $t \geq 0$*:*

$$|x(t;x_0,d)| \leq R(x_0) \tag{2.2}$$

*Moreover, suppose that for every* $\varepsilon > 0$*,* $s \geq 0$ *there exists* $T(\varepsilon, s) \geq 0$ *such that for every* $(x_0, d) \in \Re^n \times M_D$ *with* $|x_0| \leq s$ *the solution* $x(t;x_0,d)$ *of (2.1) satisfies* $|x(t;x_0,d)| \leq \varepsilon$ *for all* $t \geq T(\varepsilon, s)$ *(Uniform attractivity for bounded sets of initial states).*

*Then* $0 \in \Re^n$ *is URGAS for system (2.1).*

The following lemma provides sufficient conditions for the reachability condition (P1) of Theorem 2.2. Its proof is provided in the Appendix. Notice that we do not assume that $F(d,0)=0$ for all $d \in D$ (i.e., we do not assume the existence of an equilibrium point).

**Lemma 2.5:** *Consider system (2.1) under hypotheses (H1-3) and suppose that there exist locally Lipschitz functions $h:\Re^n \to \Re$ with $h(0)<0$, $a:\Re^n \to \Re$ being bounded from above with $a(0)=0$, $W:\Re^n \to \Re^+$ being radially unbounded, a continuous function $\delta:\Re^+ \to (0,+\infty)$ and constants $K \ge 0$, $b>0$, such that $\{x \in \Re^n : 0 < h(x) < b\} \neq \emptyset$ and*

$$\sup_{d\in D} \nabla h(x)F(d,x) \le 0 \text{, for almost all } x \in \Re^n \text{ with } 0 < h(x) < b \tag{2.3a}$$

$$\sup_{d\in D} \left(\nabla h(x) - \nabla a(x)\right)F(d,x) \le -\delta(h(x)) \text{, for almost all } x \in \Re^n \text{ with } h(x) > 0 \tag{2.3b}$$

$$\sup_{d\in D} \nabla W(x)F(d,x) \le KW(x) \text{, for almost all } x \in \Re^n \text{ with } h(x) > 0 \tag{2.4}$$

*Then for every $\hat{\varepsilon} \in (0,b)$ there exist functions $T \in C^0(\Re^n;\Re^+)$, $G \in C^0(\Re^n;\Re^+)$ such that property (P1) of Theorem 2.2 holds with $\Omega := \{x \in \Re^n : h(x) \le \hat{\varepsilon}\}$.*

We next consider the control system (1.1). The following hypotheses will be valid for system (1.1) throughout this section:

**(Q1)** $D \subset \Re^l$ is compact and $U \subseteq \Re^m$ is a convex set with $0 \in U$.

**(Q2)** The mappings $D \times \Re^n \ni (d,x) \to f(d,x) \in \Re^n$, $D \times \Re^n \ni (d,x) \to g(d,x) \in \Re^{n\times m}$ are continuous with $f(d,0)=0$ for all $d \in D$.

**(Q3)** There exists a symmetric positive definite matrix $P \in \Re^{n\times n}$ such that for every pair of compact sets $S \subset \Re^n$, $V \subseteq U$ it holds that $\sup\left\{ \frac{(x-y)'P(f(d,x)+g(d,x)u-f(d,y)-g(d,y)u)}{|x-y|^2} : d \in D, u \in V, x,y \in S, x \neq y \right\} < +\infty$.

The following theorem provides relaxed sufficient Lyapunov-like conditions for the existence of a locally Lipschitz, globally stabilizing feedback law $k:\Re^n \to U$. The Lyapunov-like conditions of the following theorem are "relaxed" in the sense that the Lyapunov differential inequality is not required to hold for every non-zero state, but only for states that belong to an appropriate set of the state space (compare with the results in [9]). On the other hand additional conditions must hold. Its proof is provided in the Appendix.

**Theorem 2.6:** *Consider system (1.1) under hypotheses (Q1-3) and suppose that there exist continuously differentiable functions $h:\Re^n \to \Re$ with $h(0)<0$, $W:\Re^n \to \Re^+$ being radially unbounded, $V:\Re^n \to \Re^+$ being positive definite and radially unbounded, a continuous non-increasing function $\delta:\Re^+ \to (0,+\infty)$ and constants $K \ge 0$, $\varepsilon > 0$ such that $\{x \in \Re^n : h(x) \ge \varepsilon\} \neq \emptyset$ and the following properties hold:*

**(R1)** *For every $x \in \Re^n$ with $h(x) \ge 0$ there exists $u \in U$ with*

$$\sup_{d\in D} \nabla h(x)\left(f(d,x)+g(d,x)u\right) \le -\delta(h(x)) \tag{2.5}$$

$$\sup_{d\in D} \nabla W(x)\left(f(d,x)+g(d,x)u\right) \le KW(x) \tag{2.6}$$

**(R2)** *For every* $x \neq 0$ *with* $h(x) \leq \varepsilon$ *there exists* $u \in U$ *with*

$$\sup_{d \in D} \nabla V(x)\big(f(d,x) + g(d,x)u\big) < 0 \tag{2.7}$$

**(R3)** *For every* $x \in \Re^n$ *with* $h(x) \in [0, \varepsilon]$ *there exists* $u \in U$ *satisfying (2.5), (2.6) and (2.7).*

**(R4)** *There exists a neighbourhood* $N$ *of* $0 \in \Re^n$ *and a locally Lipschitz mapping* $\tilde{k} : N \to U$ *with* $\tilde{k}(0) = 0$ *such that* $\sup_{d \in D} \nabla V(x)\big(f(d,x) + g(d,x)\tilde{k}(x)\big) < 0$ *for all* $x \in N$, $x \neq 0$.

*Then there exists a locally Lipschitz mapping* $k : \Re^n \to U$ *with* $k(0) = 0$ *such that* $0 \in \Re^n$ *is URGAS for the closed-loop system (1.1) with* $u = k(x)$.

**Remark 2.7:**
**(a)** It should be noted that if the mapping $\tilde{k} : N \to U$ involved in hypothesis (R4) of Theorem 2.6 is $C^1$ then the obtained feedback $k : \Re^n \to U$ is of class $C^1$ as well. Similarly, if the mappings $\tilde{k} : N \to U$ and $h : \Re^n \to \Re$ are $C^j$ ($1 \leq j \leq \infty$) then the obtained feedback $k : \Re^n \to U$ is of class $C^j$ as well.

**(b)** As already noted in the Introduction, Theorem 2.6 can be used for various purposes. For example, if $V : \Re^n \to \Re^+$ is a RCLF for (1.1), the result of Theorem 2.6 can be used in order to obtain a different family of robust feedback stabilizers from the family of robust feedback stabilizers obtained by using the classical Lyapunov feedback design methodology (see [2,9,25]). Indeed, the following section is devoted to the presentation of an important control system, for which simple formulae of robust feedback stabilizers are obtained by the use of Theorem 2.2 and Theorem 2.6, while complicated formulae of robust feedback stabilizers are obtained by the use of classical results. On the other hand, Theorem 2.2 and Theorem 2.6 can be used for the exploitation of a function $V : \Re^n \to \Re^+$ which is not necessarily a RCLF. This is the case of systems with input constraints presented in Section 4.

**(c)** Some comments concerning hypotheses (R1-4) of Theorem 2.6 are given next: hypothesis (R1) allows the construction of a feedback law which guarantees that Lemma 2.5 can be applied for the corresponding closed-loop system. Hypothesis (R2) allows the construction of a (different) feedback law which guarantees that the time derivative of the Lyapunov function is negative definite on the set $\{x \in \Re^n : h(x) < \varepsilon\}$. On the other hand, hypothesis (R3) is a crucial hypothesis that guarantees that the two feedback laws constructed by means of hypotheses (R1-2) can be combined on the region $\{x \in \Re^n : 0 < h(x) < \varepsilon\}$. Finally, hypothesis (R4) is a local hypothesis, which automatically guarantees the small-control property (see [9,25]) and allows us to construct a locally Lipschitz feedback law (instead of a simply continuous one).

## 3. Application to the Robust Global Stabilization of the Chemostat

Continuous stirred microbial bioreactors, often called chemostats, cover a wide range of applications. The dynamics of the chemostat is often adequately represented by a simple dynamic model involving two state variables, the microbial biomass $X$ and the limiting organic substrate $S$ (see [24]). For control purposes, the manipulated input is usually the dilution rate *D*. A commonly used delay-free model for microbial growth on a limiting substrate in a chemostat is of the form:

$$\begin{aligned} &\dot{X} = (\mu(S) - D)X \\ &\dot{S} = D(S_i - S) - K\mu(S)X \\ &X \in (0,+\infty), S \in (0, S_i), D \geq 0 \end{aligned} \tag{3.1}$$

where $S_i$ is the feed substrate concentration, $\mu(S)$ is the specific growth rate and $K > 0$ is a biomass yield factor. In most applications, Monod or Haldane or generalized Haldane models are used for $\mu(S)$ [3]. The reader should notice that chemostat models with time delays were considered in [8,33,34]. The literature on control studies of chemostat models of the form (3.1) is extensive. In [6], feedback control of the chemostat by manipulating the dilution rate was studied for the promotion of coexistence. Other interesting control studies of the chemostat can be found in [1,10,11,13,17,21]. The stability and robustness of periodic solutions of the chemostat was studied in [22,23]. The problem of the stabilization of a non-trivial steady state $(X_s, S_s)$ of the chemostat model (1.1) was considered in [17], where it was shown that the simple feedback law $D = \mu(S)\frac{X}{X_s}$ is a globally stabilizing feedback. See also the recent work [13] for the study of the robustness properties of the closed-loop system (3.1) with $D = \mu(S)\frac{X}{X_s}$ for time-varying inlet substrate concentration $S_i$.

In this work we consider the robust global feedback stabilization problem for the more general uncertain chemostat model

$$\begin{aligned}
&\dot{X} = \big(\mu(S) + \Delta(S,t) - D - b\big)X \\
&\dot{S} = D(S_i - S) - K\mu(S)X + m\,X \\
&X \in (0,+\infty)\,, S \in (0, S_i)\,, D \geq 0
\end{aligned} \tag{1.5}$$

In the above,

- The term $bX$ in the biomass balance represents the death rate of the cells in the chemostat. The parameter $b \geq 0$ is the cell mortality rate.
- The term $m\,X$ in the substrate balance accounts for the rate of substrate consumption for cell maintenance (see [3], pages 390 and 450) as well as the rate of release of substrate due to the death of the cells in the chemostat (which is proportional to $bX$). The parameter $m$ is either negative or assumes a small positive value. The parameter $m$ is related to the presence of variable apparent yield coefficient (which has been studied recently in [35]).
- The term $\Delta(S,t)$ represents possible deviations of the specific growth rate of the biomass, primarily accounting for the adjustment of the biomass to changes in the substrate levels. The following assumption is made about the uncertainty term $\Delta(S,t)$:

**(S0)** There exist constants $S_s \in (0, S_i)$ and $a \geq 0$ such that $\Delta(S,t) = d_1(t)|S - S_s| - d_2(t)\max\{0, S_s - S\}$, where $d_i : \Re^+ \to [0, a]$ ($i = 1,2$) are measurable, essentially bounded functions.

Clearly, $a \geq 0$ is a constant which quantifies the uncertainty range, $S_s \in (0, S_i)$ is the value of the substrate concentration where $\dot{X}$ is precisely known to be equal to $\big(\mu(S_s) - D - b\big)X$. Notice that at $S_s \in (0, S_i)$, the uncertainty $\Delta(S,t)$ is assumed to vanish. See Figure 1 below for a sketch of the shape of the uncertainty range as a function of $S$.

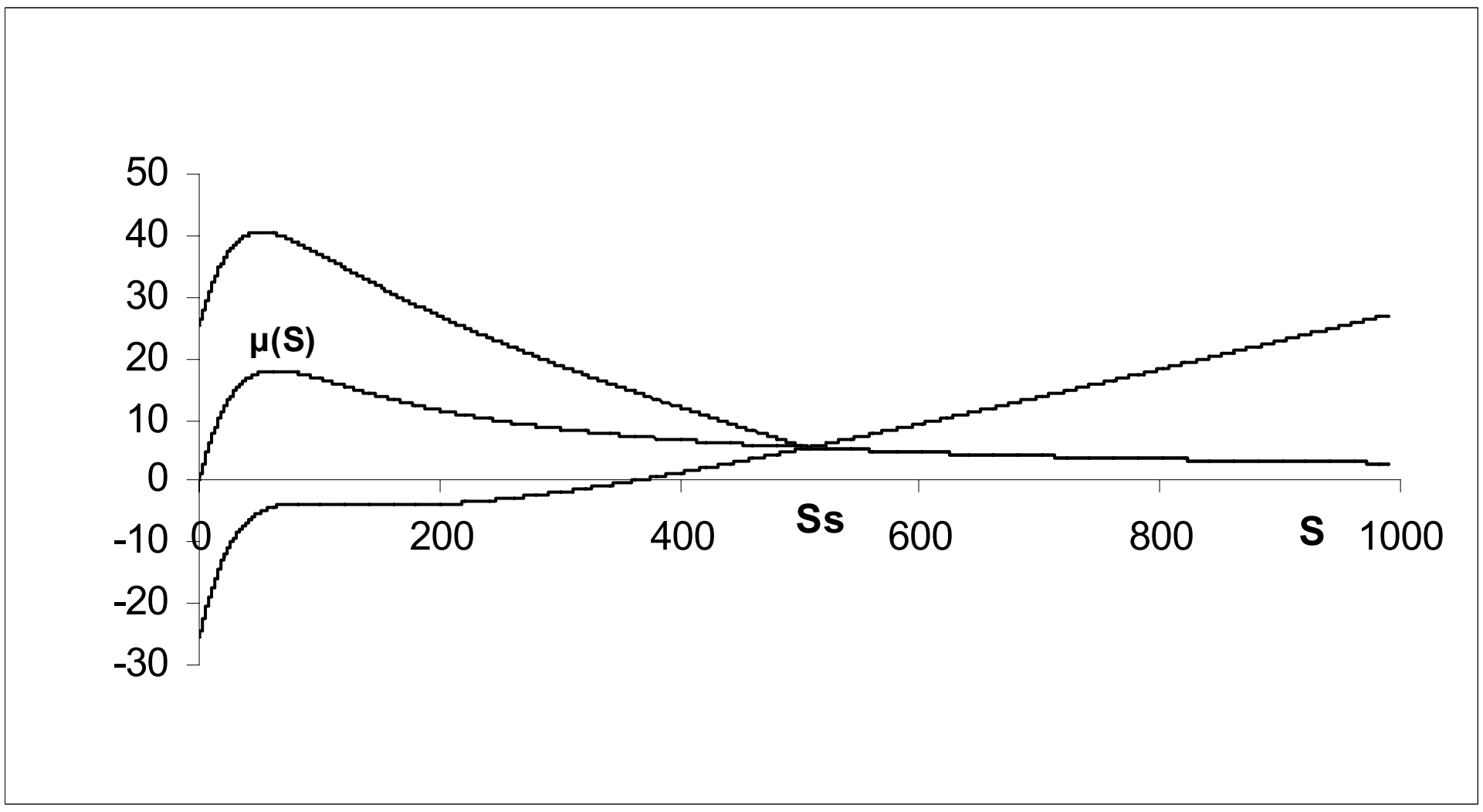

**Figure 1**: Indicative uncertainty range for the specific growth rate of the biomass $\mu(S) + \Delta(S,t)$ ( here $\mu(S) = \frac{75S}{100 + S + 0.025S^2}$, $a = 0.05$, $S_s = 506.72$ )

It should be noticed that (1.5) under hypothesis (S0) is a more general chemostat model than (3.1) (if we set $a = b = m = 0$ we obtain model (3.1)). In [10], the problem of the output regulation of the chemostat model (1.5) with $m = 0$ was considered.

It is important to notice that even in the case of zero uncertainty and zero mortality rate (i.e., $a = b = 0$) and for negative values for the constant $m$, the application of the feedback law $D = \mu(S)\frac{X}{X_s}$ does not necessarily lead to global stability. For example, for the Haldane model $\mu(S) := \frac{\mu_{\max} S}{K_1 + S + K_2 S^2}$, it is easy to verify that for arbitrarily small negative values for the constant $m$, the closed-loop system (1.5) under hypothesis (S0) with $D = \mu(S)\frac{X}{X_s}$ and $a = b = 0$ has two equilibrium points in the first quadrant with coordinates $(S_1, X_s)$ and $(S_2, X_s)$, where $0 < S_1 < S_2$. The equilibrium point $(S_2, X_s)$ is locally asymptotically stable with region of attraction the set $\{(S, X) : S > S_1, X > 0\}$. The stable manifold of the unstable equilibrium $(S_1, X_s)$ is the straight line $S = S_1$ and if the initial condition for the substrate is less than $S_1$ then the system is led to shut-down in finite time (i.e., there exists $T \geq 0$ such that $\lim_{t \to T^-} S(t) = 0$). Therefore, the feedback law $D = \mu(S)\frac{X}{X_s}$ needs to be modified in order to be able to guarantee global asymptotic stability for the desired equilibrium point.

Throughout this section we will assume that the specific growth rate function $\mu : \Re \to [0, \mu_{\max}]$ involved in the chemostat models (1.5) is a locally Lipschitz function with $\mu(S) = 0$ for all $S \leq 0$ and $\mu(S) > 0$ for all $S > 0$. We consider system (1.5) under the following additional hypotheses:

**(S1)** There exists an equilibrium point $(X_s, S_s) \in (0, +\infty) \times (0, S_i)$ with $\mu(S_s) = D_s + b$ and $\frac{D_s (S_i - S_s)}{K(D_s + b) - m} = X_s$ for certain value of the dilution rate $D_s > 0$.

Assumption (S1) is satisfied for Monod, Haldane and generalized Haldane kinetics, as long as the value of the dilution rate $D_s$ is not too high.

**(S2)** There exists $S^+ \in (0, S_s)$ and $p > 0$ such that $K\mu(S) - m \geq p$ and $\mu(S) - b \geq 2p$ for all $S \in [S^+, S_i]$.

Assumption (S2) is satisfied for Monod, Haldane and generalized Haldane kinetics, as long as $\min\left(\mu(S_s), \mu(S_i)\right) > \max\left(b, \frac{m}{K}\right)$.

The goal is the robust global stabilization of the non-trivial equilibrium point $(X_s, S_s) \in (0, +\infty) \times (0, S_i)$ with $\mu(S_s) = D_s + b$ and $\frac{D_s (S_i - S_s)}{K(D_s + b) - m} = X_s$ involved in hypotheses (S1-2) for system (1.5). To this end we apply the change of coordinates:

$$S = \frac{S_i \exp(x_1)}{c + \exp(x_1)} \; ; \; \frac{X}{S_i - S} = G \exp(x_2) \tag{3.2}$$

and the input transformation:

$$D = D_s + u \tag{3.3}$$

where $c:=\frac{S_i}{S_s}-1$ and $G:=\frac{D_s}{K(D_s+b)-m}$. The above coordinate change maps the strip $\left\{(X,S)\in\Re^2 : X>0\,,0<S<S_i\right\}$ onto $\Re^2$. Under the above transformation system (1.5) under hypothesis (S0) is expressed by the following control system:

$$
\begin{aligned}
\dot{x}_1 &= \left(c\exp(-x_1)+1\right)\left(D_s+u-(K\tilde{\mu}(x_1)-m)G\exp(x_2)\right) \\
\dot{x}_2 &= \left(\tilde{\mu}(x_1)+d_1\frac{cS_s}{c+\exp(x_1)}\left|\exp(x_1)-1\right|-d_2\frac{cS_s}{c+\exp(x_1)}\max\left(0,1-\exp(x_1)\right)-b\right)-(K\tilde{\mu}(x_1)-m)G\exp(x_2) \\
x &= (x_1,x_2)\in\Re^2, u\in U:=[-D_s,+\infty), d=(d_1,d_2)\in[0,a]^2
\end{aligned}
\tag{3.4}
$$

where $\tilde{\mu}(x_1):=\mu\left(\frac{S_i\exp(x_1)}{c+\exp(x_1)}\right)$. Notice that $x=0$ is an equilibrium point for the above system for $u\equiv 0$. Therefore, we seek for a locally Lipschitz feedback law $k:\Re^2\to\Re$ with $k(0)=0$ so that $0\in\Re^2$ is uniformly robustly globally asymptotically stable for the closed-loop system (3.4) with $u=k(x)$ in the sense described in the previous section.

Insights for the solution of the feedback stabilization problem for (3.4) may be obtained by setting $a=b=m=0$ and obtaining the transformed system (3.1):

$$
\begin{aligned}
\dot{x}_1 &= \left(c\exp(-x_1)+1\right)\left(D_s+u-\tilde{\mu}(x_1)\exp(x_2)\right) \\
\dot{x}_2 &= \tilde{\mu}(x_1)-\tilde{\mu}(x_1)\exp(x_2) \\
x &= (x_1,x_2)\in\Re^2, u\in U:=[-D_s,+\infty)
\end{aligned}
\tag{3.5}
$$

For the control system (3.5), families of CLFs are known (see [13]). Let $\gamma:\Re\to\Re^+$, $\beta:\Re\to\Re^+$ be non-negative, continuously differentiable functions with $\gamma(0)=\beta(0)=0$ and such that

$$x\gamma'(x)>0\,,\ x\beta'(x)>0\,,\ \text{for all } x\neq 0 \tag{3.6a}$$

$$\text{if } x\to\pm\infty \text{ then } \gamma(x)\to+\infty \text{ and } \beta(x)\to+\infty \tag{3.6b}$$

For example the functions $\gamma(x)$ and $\beta(x)$ could be of the form $Kx^m$, where $K>0$ and $m>0$ is an even positive integer. Inequalities (3.6a,b) guarantee that the following family of functions:

$$V(x)=\gamma(x_1)+\beta(x_2) \tag{3.7}$$

are radially unbounded, positive definite and continuously differentiable functions. The reader may verify that the above functions are CLFs for the control system (3.5). The knowledge of the above family of CLFs allows us to obtain a family of stabilizing feedback laws for (3.5). The reader may verify that the following family of feedback laws:

$$k(x):=-D_s+\tilde{\mu}(x_1)\exp(x_2)\varphi(x_1)+q(x_1,x_2) \tag{3.8}$$

where $\varphi:\Re\to\Re^+$ is a locally Lipschitz, non-negative function with $\varphi(0)=1$, $\varphi(x)<1$ for $x>0$ and $\varphi(x)>1$ for $x<0$, $q:\Re^2\to\Re^+$ is a locally Lipschitz, non-negative function with $q(x_1,x_2)=0$ for $x_1\geq 0$, is a family of globally stabilizing feedback laws for (3.5). For example, the selection $\varphi(x)=\frac{c+1}{c+\exp(x)}$, $q\equiv 0$ gives a feedback law, which transformed back to the original coordinates gives $D=\mu(S)\frac{X}{X_s}$. This is the feedback law considered in [17] (see also [13] and references therein).

On the other hand, it can be verified that $V$ as defined by (3.7), where $\gamma:\Re\to\Re^+$ and $\beta:\Re\to\Re^+$ satisfy (3.6), is not necessarily a CLF for (3.4) under hypotheses (S1-2). However, we will show next that $\gamma:\Re\to\Re^+$ and $\beta:\Re\to\Re^+$ can be selected so that $V$ as defined by (3.7) is a CLF for (3.4) under hypotheses (S1-2). Hypothesis (S2) implies that there exists $x_1^* \in [x_1^+,0)$, where $x_1^+ = \ln\left(\frac{S^+ c}{S_i - S^+}\right)$, such that

$$\frac{acS_s}{c+\exp(x_1)}\max\left(0,1-\exp(x_1)\right)\le p\ ,\ K\tilde{\mu}(x_1)-m\ge p \text{ and } \tilde{\mu}(x_1)-b\ge 2p \text{ for all } x_1\ge x_1^* \tag{3.9}$$

**Proposition 3.1:** *The function*

$$V(x):=\gamma(x_1)+\frac{1}{2}x_2^2 \tag{3.10}$$

*where*

$$\gamma(x_1):=\frac{1}{2}M\,x_1^2+\begin{cases} A\left[\exp(2(x_1+x_1^*))-1-2(x_1+x_1^*)\right] & ,\quad for\ x_1\ge -x_1^* \\ 0 & ,\quad for\ x_1< -x_1^* \end{cases} \tag{3.11}$$

*and $M, A>0$ are constants sufficiently large, is a RCLF for (3.4) under hypotheses (S1-2). Moreover, the feedback law:*

$$k(x)=-D_s+\max\{0\,,K\tilde{\mu}(x_1)-m\}\,G\exp(x_2-x_1)+q(x_1,x_2) \tag{3.12}$$

*where $q:\Re^2\to\Re^+$ is a locally Lipschitz, non-negative function with $q(x_1,x_2)=0$ for $x_1\ge 0$ and*

$$q(x_1,x_2)\ge\frac{W(x_1,x_2)\exp(x_1)}{M\left(c+\exp(x_1)\right)}+a\frac{cS_s|x_2|\exp(x_1)}{M(c+\exp(x_1))^2}-\exp(x_1)x_2\frac{\tilde{\mu}(x_1)-b-(K\tilde{\mu}(x_1)-m)G\exp(x_2)}{Mx_1\left(c+\exp(x_1)\right)}\ ,\text{ for } x_1\le x_1^* \tag{3.13}$$

*where $W:\Re^2\to\Re^+$ is a locally Lipschitz, positive definite function, satisfies for all $x\neq 0$, $d=(d_1,d_2)\in[0,a]^2$:*

$$\dot V=\nabla V(x)\begin{bmatrix} \left(c\exp(-x_1)+1\right)\left(D_s+k(x)-(K\tilde{\mu}(x_1)-m)G\exp(x_2)\right) \\ \tilde{\mu}(x_1)+d_1\frac{cS_s}{c+\exp(x_1)}\left|\exp(x_1)-1\right|-d_2\frac{cS_s}{c+\exp(x_1)}\max\left(0,1-\exp(x_1)\right)-b-(K\tilde{\mu}(x_1)-m)G\exp(x_2)\end{bmatrix}<0$$

**Proof:** Notice that the function $V$ defined by (3.10), (3.11) is continuously differentiable, positive definite and radially unbounded.

Let $L:=\sup\left\{\frac{|\mu(S)-\mu(S_s)|}{|S-S_s|}:S\in(0,S_i),S\neq S_s\right\}$. Clearly, $L\le L_\mu<+\infty$ where $L_\mu$ denotes the Lipschitz constant for the specific growth rate function on $[0,S_i]$. Using definition $\tilde{\mu}(x_1):=\mu\left(\frac{S_i\exp(x_1)}{c+\exp(x_1)}\right)$, we obtain:

$$\left|\tilde{\mu}(x_1)-\tilde{\mu}(0)\right|\le LS_s c\frac{\left|\exp(x_1)-1\right|}{c+\exp(x_1)}\ ,\text{ for all } x_1\in\Re \tag{3.14}$$

Let $\delta\in(0,p)$ and define:

$$\beta_{\min} := \ln\left(\frac{p-\delta}{(K\mu_{\max}-m)G}\right);\ \beta_{\max} := \ln\left(\frac{1}{pG}\left(\mu_{\max}+a(c+1)S_s-b+\delta\right)\right) \tag{3.15}$$

Notice that by selecting $\delta \in (0,p)$ sufficiently close to $p$ we have $\beta_{\min} < 0 < \beta_{\max}$. We will show next that the function $V$ defined by (3.10), (3.11) is a RCLF for (3.4) for constants $M, A > 0$ that satisfy:

$$2AcpG\exp(\beta_{\min})\exp(2x_1^*) \geq BS_s\left(L\frac{|Kb-m|}{K(D_s+b)-m}+a\right)+1 \tag{3.16a}$$

$$M \geq \frac{2A}{-x_1^*};\ M \geq \frac{1}{2c}\exp(-\beta_{\min}-x_1^*)+\frac{c}{2}\exp(-\beta_{\min}-x_1^*)\left(\frac{S_s r}{\varepsilon\, pG}\right)^2\left(\frac{L|Kb-m|}{K(D_s+b)-m}+2a\right)^2 \tag{3.16b}$$

where $B := \max\left\{|\beta_{\min}|, |\beta_{\max}|\right\}$ and $\varepsilon > 0$ is sufficiently small such that

$$x_1(\exp(-x_1)-1) \leq -\varepsilon x_1^2 \text{, for all } x_1 \in [x_1^*, -x_1^*] \tag{3.17a}$$

$$x_2(1-\exp(x_2)) \leq -\varepsilon x_2^2 \text{, for all } x_2 \in [\beta_{\min}, \beta_{\max}] \tag{3.17b}$$

We consider the following cases.

**Case 1:** $x_1 \geq x_1^*$ and $x_2 \notin [\beta_{\min}, \beta_{\max}]$.

Notice that by virtue of (3.9) and definitions (3.15), the following inequalities hold for $x_1 \geq x_1^*$:

$$\beta_{\min} \leq \ln\left(\frac{1}{(K\tilde{\mu}(x_1)-m)G}\left(\tilde{\mu}(x_1)-b-\frac{acS_s}{c+\exp(x_1)}\max\left(0,1-\exp(x_1)\right)-\delta\right)\right) \tag{3.18a}$$

$$\beta_{\max} \geq \ln\left(\frac{1}{(K\tilde{\mu}(x_1)-m)G}\left(\tilde{\mu}(x_1)+\frac{acS_s}{c+\exp(x_1)}|\exp(x_1)-1|-b+\delta\right)\right) \tag{3.18b}$$

Using (3.18a), we obtain for all $x_1 \geq x_1^*$, $x_2 \leq \beta_{\min}$, $d=(d_1,d_2) \in [0,a]^2$:

$$\left(\tilde{\mu}(x_1)+d_1\frac{cS_s}{c+\exp(x_1)}|\exp(x_1)-1|-d_2\frac{cS_s}{c+\exp(x_1)}\max\left(0,1-\exp(x_1)\right)-b\right)-(K\tilde{\mu}(x_1)-m)G\exp(x_2) \geq \delta \tag{3.19a}$$

Using (3.18b), we obtain for all $x_1 \geq x_1^*$, $x_2 \geq \beta_{\max}$, $d=(d_1,d_2) \in [0,a]^2$:

$$\left(\tilde{\mu}(x_1)+d_1\frac{cS_s}{c+\exp(x_1)}|\exp(x_1)-1|-d_2\frac{cS_s}{c+\exp(x_1)}\max\left(0,1-\exp(x_1)\right)-b\right)-(K\tilde{\mu}(x_1)-m)G\exp(x_2) \leq -\delta \tag{3.19b}$$

Consequently, we obtain from (3.12) and (3.19a,b) for all $x_1 \geq x_1^*$, $x_2 \notin [\beta_{\min}, \beta_{\max}]$ and $d=(d_1,d_2) \in [0,a]^2$:

$$\dot{V} \leq \gamma'(x_1)\left(c\exp(-x_1)+1\right)pG\exp(x_2)\left(\exp(-x_1)-1\right)-\delta|x_2| < 0 \tag{3.20}$$

**Case 2:** $x_1 \in [x_1^*, -x_1^*]$ and $x_2 \in [\beta_{\min}, \beta_{\max}]$.

In this case, using (3.12), we have:

$$\dot{V} \leq Mx_1\left(c\exp(-x_1)+1\right)\left(K\tilde{\mu}(x_1)-m\right)G\exp(x_2)\left(\exp(-x_1)-1\right)$$
$$+\frac{Kb-m}{K(D_s+b)-m}x_2\left(\tilde{\mu}(x_1)-\tilde{\mu}(0)\right)+\frac{cS_s}{c+\exp(x_1)}x_2 d_1\left|\exp(x_1)-1\right|-d_2\frac{cS_s}{c+\exp(x_1)}x_2\max\left(0,1-\exp(x_1)\right) \quad (3.21)$$
$$+(K\tilde{\mu}(x_1)-m)G\,x_2(1-\exp(x_2))$$

Let $r>0$ such that $\frac{\left|\exp(x_1)-1\right|}{c+\exp(x_1)} \leq r|x_1|$ for all $x_1 \in [x_1^*,-x_1^*]$. Using the previous inequality in conjunction with (3.9), (3.14), (3.17a,b) and (3.21) we obtain for all $x_1 \in [x_1^*,-x_1^*]$, $x_2 \in [\beta_{\min},\beta_{\max}]$ and $d=(d_1,d_2)\in[0,a]^2$ :

$$\dot{V} \leq -\varepsilon M\left(c+\exp(x_1^*)\right)pG\exp(\beta_{\min}+x_1^*)\,x_1^2+\left(\frac{L|Kb-m|}{K(D_s+b)-m}+2a\right)cS_s r\,|x_2||x_1|-\varepsilon\,pG\,x_2^2 \quad (3.22)$$

By completing the squares in the right-hand side of (3.22) and using (3.16b), we obtain for all $x_1 \in [x_1^*,-x_1^*]$, $x_2 \in [\beta_{\min},\beta_{\max}]$ and $d=(d_1,d_2)\in[0,a]^2$ :

$$\dot{V} \leq -\frac{\varepsilon\,pG}{2}\left(x_1^2+x_2^2\right) \quad (3.23)$$

**Case 3:** $x_1 \geq -x_1^*$ and $x_2 \in [\beta_{\min},\beta_{\max}]$.

Let $B:=\max\left\{|\beta_{\min}|,|\beta_{\max}|\right\}$. In this case, using (3.9), (3.12) and (3.14) we obtain for all $x_1 \geq -x_1^*$, $x_2 \in [\beta_{\min},\beta_{\max}]$ and $d=(d_1,d_2)\in[0,a]^2$ :

$$\dot{V} \leq \left[\gamma'(x_1)cpG\exp(\beta_{\min}-2x_1)-BS_s\left(L\frac{|Kb-m|}{K(D_s+b)-m}+a\right)\right]\left(1-\exp(x_1)\right)+pG\,x_2(1-\exp(x_2)) \quad (3.24)$$

By virtue of definition (3.11) we get for all $x_1 \geq -x_1^*$ :

$$\gamma'(x_1)cpG\exp(\beta_{\min}-2x_1)=2AcpG\exp(\beta_{\min})\exp(2x_1^*)+cpG\exp(\beta_{\min}-2x_1)\left[M\,x_1-2A\right] \quad (3.25)$$

Using (3.16a,b) and (3.25), we obtain for all $x_1 \geq -x_1^*$, $x_2 \in [\beta_{\min},\beta_{\max}]$ and $d=(d_1,d_2)\in[0,a]^2$ :

$$\dot{V} \leq \left(1-\exp(x_1)\right)+pG\,x_2(1-\exp(x_2))<0 \quad (3.26)$$

**Case 4:** $x_1 \leq x_1^*$

In this case we have:

$$\dot{V}=Mx_1\left(c\exp(-x_1)+1\right)\left(D_s+k(x)-(K\tilde{\mu}(x_1)-m)G\exp(x_2)\right)$$
$$+x_2\left(\tilde{\mu}(x_1)+(d_1-d_2)\frac{cS_s}{c+\exp(x_1)}\left(1-\exp(x_1)\right)-b-(K\tilde{\mu}(x_1)-m)G\exp(x_2)\right) \quad (3.27)$$

By virtue of (3.12), (3.13) and (3.27) we obtain for all $x_1 \leq x_1^*$, $x_2 \in \Re$ and $d=(d_1,d_2)\in[0,a]^2$ :

$$\dot{V} \leq x_1 W(x_1,x_2)<0 \quad (3.28)$$

The proof is complete. ⊲

We next consider the possibility of constructing simpler feedback laws than the family of feedback laws given by (3.12), (3.13). To this purpose we utilize the relaxed Lyapunov-like conditions of Theorem 2.6 and the stability conditions of Theorem 2.2.

**Theorem 3.2:** *Let $\psi:\Re\to\Re^+$ be a locally Lipschitz non-increasing function with $\psi(s)=0$ for all $s\geq 0$ and $\psi(s)>0$ for all $s<0$ and let $L:\Re^2\to(0,+\infty)$ be a locally Lipschitz function with $\inf\{L(x):x\in\Re^2\}>0$. Under hypotheses (S1-2), for every $a\geq 0$, $0\in\Re^2$ is URGAS for the closed-loop system (3.4) with*

$$u=-D_s+\max(0,K\tilde{\mu}(x_1)-m)G\exp(x_2-x_1)+L(x_1,x_2)\psi(x_1) \tag{3.29}$$

**Proof:** We define:

$$\Omega:=\{(x_1,x_2)\in\Re^2:x_1\geq x_1^*\} \tag{3.30}$$

Following exactly the same arguments as in cases 1,2,3 of the proof of Proposition 3.1, we are in a position to show that hypothesis (P2) of Theorem 2.2 holds for the closed-loop system (3.4) with (3.29) with $\Omega$ as defined by (3.30) and $V$ as defined by (3.10), (3.11) for sufficiently large constants $M,A>0$. Next define:

$$h(x):=\frac{1}{2}x_1^*-x_1 \tag{3.31}$$

It should be noticed that by virtue of definitions (3.30) and (3.31), the set $\Omega$ satisfies $\Omega:=\{x\in\Re^2:h(x)\leq\hat{\varepsilon}\}$ with $\hat{\varepsilon}=-\frac{1}{2}x_1^*>0$. Let $l:=\inf\{L(x):x\in\Re^2\}>0$. It follows from (3.4), (3.29), (3.31) and the fact that $\psi:\Re\to\Re^+$ is non-increasing, that the following inequality holds for all $x\in\Re^2$ with $h(x)\geq 0$:

$$\begin{aligned}\dot{h}=\nabla h(x)\dot{x}&=-(c\exp(-x_1)+1)(D_s+u-(K\tilde{\mu}(x_1)-m)G\exp(x_2))\\&\leq -l\left(c\exp(-\frac{1}{2}x_1^*)+1\right)\left(1-\psi\left(\frac{1}{2}x_1^*\right)\right)\end{aligned} \tag{3.32}$$

Consequently, inequalities (2.3a,b) of Lemma 2.5 holds for every $b>\hat{\varepsilon}$ with $a(x)\equiv 0$ and

$$\delta(h(x)):=l\left(c\exp(-\frac{1}{2}x_1^*)+1\right)\left(1-\psi\left(\frac{1}{2}x_1^*\right)\right)>0 \tag{3.33}$$

Finally, define the continuously differentiable function:

$$W(x_1,x_2):=\frac{1}{2}x_1^2+\frac{1}{2}(\min\{0,x_2\})^2+\frac{1}{2}\left(\max\{0,x_2-\ln(c+e^{x_1})\}\right)^2+1 \tag{3.34}$$

Notice that the following inequalities hold for all $x_1,x_2\leq 0$ and $d=(d_1,d_2)\in[0,a]^2$:

$$-aS_s-K\mu_{\max}G\leq\tilde{\mu}(x_1)+(d_1-d_2)\frac{cS_s}{c+\exp(x_1)}(1-\exp(x_1))-(K\tilde{\mu}(x_1)-m)G\exp(x_2)\leq\mu_{\max}+aS_s+mG \tag{3.35}$$

as well as the following equality:

$$\frac{d}{dt}\left(x_2-\ln(c+e^{x_1})\right)=\tilde{\mu}(x_1)+d_1\frac{cS_s}{c+\exp(x_1)}|\exp(x_1)-1|-d_2\frac{cS_s}{c+\exp(x_1)}\max(0,1-\exp(x_1))-b-D_s-u \tag{3.36}$$

Using inequalities (3.35) in conjunction with (3.34), (3.35) we obtain inequality (2.4) for certain constant $K>0$ sufficiently large. Consequently, Lemma 2.5 implies that hypothesis (P1) of Theorem 2.2 holds for the closed-loop system (3.4) with (3.29) with $\Omega$ as defined by (3.30). Theorem 2.2 implies that $0 \in \Re^2$ is URGAS for the closed-loop system (3.4) with (3.29). The proof is complete. $\lhd$

It should be noticed that the obtained family of stabilizing feedback laws (3.29) is much simpler than the family of stabilizing feedback laws given by (3.12), (3.13). Moreover, it should be emphasized that the family of stabilizing feedback laws (3.29) and the family of stabilizing feedback laws given by (3.12), (3.13) do not coincide (although both families of feedback laws are members of the family expressed by (3.8) for the case $m = 0$). The reader should notice that the feedback law (3.29) transformed back to the original coordinates is expressed by (1.6). Finally, it is clear that the feedback law (3.29) is independent of the constant $a \geq 0$ which quantifies the uncertainty range. Therefore, the feedback law (3.29) achieves stabilization of $0 \in \Re^2$ for all $a \geq 0$, i.e., for arbitrary large range of uncertainty.

## 4. Feedback Stabilization of Control Systems with Input Restrictions

In this section some examples are provided, which show that the notion of relaxed Control Lyapunov Functions is very useful when trying to design stabilizing feedback laws for control systems with input restrictions. Our first example deals with a single input affine control system.

**Example 4.1:** Consider the single input affine in the control disturbance free system:

$$\begin{aligned} &\dot{x} = f(x) + g(x)u \\ &x \in \Re^n , u \in U := [-a, +\infty) \subset \Re \end{aligned} \tag{4.1}$$

where $a > 0$ is a constant and $f, g : \Re^n \to \Re^n$ are smooth vector fields with $f(0) = 0$. Suppose that a smooth, positive definite and radially unbounded function $V : \Re^n \to \Re^+$ is known such that

$$\nabla V(x) f(x) - \gamma(x) \big( \nabla V(x) g(x) \big)^2 < 0 \,, \ \forall x \neq 0 \tag{4.2}$$

where $\gamma : \Re^n \to \Re^+$ is a smooth function. Notice that under hypothesis (4.2), it follows that if the control input $u$ were allowed to take values in $\Re$ then $V : \Re^n \to \Re^+$ would be a CLF for (4.1) and a smooth stabilizing feedback for (4.1) would be $k(x) := -\gamma(x) \nabla V(x) g(x)$. On the other hand, the control input $u$ is restricted to take values in $U := [-a, +\infty) \subset \Re$. Clearly, the use of the feedback $k(x) := -\gamma(x) \nabla V(x) g(x)$ becomes problematic on the set $\{ x \in \Re^n : \gamma(x) \nabla V(x) g(x) > a \}$ and $V : \Re^n \to \Re^+$ is not necessarily a CLF for (4.1).

Let $\varepsilon \in (0, a)$ and define $h(x) := \gamma(x) \nabla V(x) g(x) - a + \varepsilon$. Clearly, $h(0) < 0$ and by virtue of (4.2) the smooth feedback $k(x) := -\gamma(x) \nabla V(x) g(x)$ can be applied on the set $\Omega := \{ x \in \Re^n : h(x) \leq \varepsilon \}$, i.e. hypotheses (R2), (R4) of Theorem 2.6 hold. Moreover, assume the existence of a function $W : \Re^n \to \Re^+$ being radially unbounded and a constant $K \geq 0$ such that

$$\nabla W(x) f(x) + u \nabla W(x) g(x) \leq K W(x) \text{, for all } u \leq -a + \varepsilon \text{ and } x \in \Re^n \text{ with } a - \varepsilon \leq \gamma(x) \nabla V(x) g(x) \leq a \tag{4.3}$$

$$\nabla W(x) f(x) - a \nabla W(x) g(x) \leq K W(x) \text{, for all } x \in \Re^n \text{ with } a \leq \gamma(x) \nabla V(x) g(x) \tag{4.4}$$

Furthermore, assume the existence of a positive constant $\delta > 0$ such that

$$\nabla h(x) f(x) - a \nabla h(x) g(x) \leq -\delta \text{, for all } x \in \Re^n \text{ with } a \leq \gamma(x) \nabla V(x) g(x) \tag{4.5}$$

$$\nabla h(x) f(x) + u \nabla h(x) g(x) \leq -\delta \text{, for all } u \leq -a + \varepsilon \text{ and } x \in \Re^n \text{ with } a - \varepsilon \leq \gamma(x) \nabla V(x) g(x) \leq a \tag{4.6}$$

Notice that inequalities (4.3), (4.4), (4.5), (4.6) in conjunction with inequality (4.2) guarantee that hypotheses (R1), (R3) of Theorem 2.6 hold as well. In this case an explicit formula for a locally Lipschitz feedback stabilizer can be given. The locally Lipschitz feedback law:

$$\widetilde{k}(x) := -\min\{ a ; \gamma(x) \nabla V(x) g(x) \} \tag{4.7}$$

guarantees global stabilization of $0 \in \Re^n$ for system (4.1). This fact follows from Lemma 2.5 and Theorem 2.2 in conjunction with inequalities (4.2), (4.3), (4.4), (4.5) and (4.6).

For example, the linear planar system

$$\begin{aligned} &\dot{x}_1 = -x_1 + x_2 \\ &\dot{x}_2 = u \\ &x = (x_1, x_2)' \in \Re^2 , u \in [-1, +\infty) \end{aligned} \tag{4.8}$$

satisfies all the above requirements with $V(x) := \frac{1}{2} x_1^2 + \frac{1}{2} x_2^2$, $W(x) \equiv V(x)$, $a := -1$, $\varepsilon := \frac{1}{2}$, $K := 1$, $\delta := \frac{1}{2}$ and $\gamma(x) \equiv 1$. Notice that $V(x) := \frac{1}{2} x_1^2 + \frac{1}{2} x_2^2$ is not a CLF for (4.8). It follows that the feedback law $\widetilde{k}(x) := -\min\{1 ; x_2\}$ guarantees global stabilization of $0 \in \Re^2$ for system (4.8). ⊲

The following example deals with the design of bounded feedback stabilizers for nonlinear uncertain systems. Many researchers have studied the problem of existence and design of robust bounded feedback stabilizers for control systems (see [18,19,20,28,29,30,32]). Here, we show that the "relaxed" Lyapunov conditions given in the present work can be used in order to rediscover sufficient conditions for the existence of robust bounded feedback stabilizers which have been obtained previously (see [32]).

**Example 4.2:** Here we study the problem of "adding an integrator" with bounded feedback. Particularly, we consider the system:

$$\begin{aligned} &\dot{x} = F(d, x, y) \\ &x \in \Re^n , y \in \Re , d \in D \end{aligned} \tag{4.9}$$

where $D \subset \Re^k$ is a compact set , $F : D \times \Re^n \times \Re \to \Re^n$ is a locally Lipschitz mapping with $F(d,0,0) = 0$ for all $d \in D$. We will assume next that system (4.9) can be stabilized by a locally Lipschitz bounded feedback law $y = \varphi(x)$. However, we will not assume the knowledge of a Lyapunov function for the closed-loop system (4.9) with $y = \varphi(x)$: instead we will assume the knowledge of a quadratic "relaxed" Lyapunov function for the closed-loop system (4.9) with $y = \varphi(x)$. The following set of assumptions is similar to the one presented in [32]:

**(W1)** *There exists a symmetric, positive definite matrix $P \in \Re^{n \times n}$, constants $\mu, a > 0$, a locally Lipschitz function $\varphi : \Re^n \to [-a, a]$, a vector $k \in \Re^n$ and a compact set $S \subset \Re^n$ containing a neighborhood of $0 \in \Re^n$ such that*:

$$x' P F(d, x, k'x) \le -\mu \, x' P x , \ \forall x \in S \tag{4.10}$$

$$\varphi(x) = k'x , \ \forall x \in S \tag{4.11}$$

**(W2)** *There exists a constant $c > 0$ and continuous mappings $T, Q : \Re^n \to \Re^+$ such that for all $(x_0, d, v) \in \Re^n \times M_D \times M_{[-c,c]}$ there exists $\hat{t}(x_0, d, v) \in [0, T(x_0)]$ with the property that the solution $x(t)$ of (4.9) with $y = \varphi(x) + v$, $x(0) = x_0$ corresponding to inputs $(d, v) \in M_D \times M_{[-c,c]}$ exists for all $t \ge 0$ and satisfies $x(t) \in S$ for all $t \ge \hat{t}(x_0, d, v)$, $|x(t)| \le Q(x_0)$ for all $t \in [0, \hat{t}(x_0, d, v)]$.*

**(W3)** *There exist constants* $C, b \ge 0$ *such that* $|F(d,x,\varphi(x))| \le C(1+|x|)$, $|F(d,x,\varphi(x)+v) - F(d,x,\varphi(x))| \le C|v|$, *for all* $(d,x,v) \in D \times \Re^n \times \Re$. *Moreover, it holds that* $|\nabla\varphi(x) F(d,x,\varphi(x)+v)| \le b$, *for almost all* $x \in \Re^n$ *and all* $(d,v) \in D \times [-c,c]$.

The reader should notice that by virtue of hypothesis (W1) it follows that property (P2) of Theorem 2.2 holds with $V(x) = x'Px$. Moreover hypothesis (W2) guarantees that property (P1) of Theorem 2.2 holds as well for the closed-loop system (4.9) with $y = \varphi(x)$. Therefore, Theorem 2.2 implies that $0 \in \Re^n$ is RGAS for the closed-loop system (4.9) with $y = \varphi(x)$ under hypotheses (W1-2).

Next we consider the subsystem:

$$\begin{aligned} &\dot{y} = f(d,x,y) + g(d,x,y)\gamma(u) \\ &u \in \Re \end{aligned} \tag{4.12}$$

where $f : D \times \Re^n \times \Re \to \Re^n$, $g : D \times \Re^n \times \Re \to \Re^n$, $\gamma : \Re \to \Re$ are locally Lipschitz mappings with $f(d,0,0) = 0$ for all $d \in D$, which satisfy the following hypothesis:

**(W4)** *There exist constants* $q, L, r > 0$ *such that* $|f(d,x,y)| \le \min\{q, L|x| + L|y|\}$, $r \le g(d,x,y)$, *for all* $(d,x,y) \in D \times \Re^n \times \Re$. *The mapping* $\gamma : \Re \to \Re$ *is non-decreasing, globally Lipschitz with unit Lipschitz constant and satisfies* $\gamma(u) = u$ *for all* $u \in [-\eta, \eta]$, *where* $\eta > r^{-1}(q+b)$, *where* $b \ge 0$ *is the constant involved in (W3).*

Exploiting the results of Section 2, we are in a position to prove the following lemma. Its proof is provided at the Appendix. It should be emphasized that the proof of Lemma 4.3 is based on the result of Lemma 2.5.

**Lemma 4.3:** *Consider system (4.9), (4.12) under hypotheses (W1-4). For every* $\tilde{c} \in [0, \eta - r^{-1}(q+b))$, $\tilde{a} \in (\tilde{c} + r^{-1}(q+b), \eta]$ *there exists* $p > 0$ *sufficiently large and continuous mappings* $\tilde{T}, \tilde{Q} : \Re^n \times \Re \to \Re^+$ *such that hypotheses (W1-2) hold with* $\tilde{c} \in [0, \eta - r^{-1}(q+b))$, $\tilde{a} \in (\tilde{c} + r^{-1}(q+b), \eta]$, $\tilde{x} := (x,y) \in \Re^{n+1}$, $\tilde{y} := u$, $\tilde{\varphi}(\tilde{x}) := -\tilde{a}\, sat(p(y - \varphi(x)))$, $\tilde{S} := \{(x,y) \in S \times \Re : |y - \varphi(x)| \le p^{-1}\} \subset \Re^{n+1}$, $\tilde{k} := -\tilde{a}p(-k', 1)' \in \Re^{n+1}$, $\tilde{P} := \begin{bmatrix} P + kk' & -k \\ -k' & 1 \end{bmatrix} \in \Re^{(n+1)\times(n+1)}$, $\tilde{F}(d,\tilde{x},\tilde{y}) := \begin{bmatrix} F(d,x,y) \\ f(d,x,y) + g(d,x,y)u \end{bmatrix}$ *and* $\tilde{\mu} := \frac{1}{2}\mu$, *in place of* $c > 0$, $a > 0$, $x \in \Re^n$, $y \in \Re$, $\varphi(x)$, $S \subset \Re^n$, $k \in \Re^n$, $P \in \Re^{n \times n}$, $F(d,x,y)$ *and* $\mu$, *respectively.*

*Moreover, if there exists constant* $R > 0$ *such that* $g(d,x,y) \le R$, *for all* $(d,x,y) \in D \times \Re^n \times \Re$ *then hypothesis (W3) holds as well with* $\tilde{C} := 2(C+L) + C(a+1) + R(\tilde{a}+1)$, $\tilde{b} := \tilde{a}p(q + R\tilde{a} + R\tilde{c} + b)$ *in place of* $C, b \ge 0$.

Applying induction and the result of Lemma 4.3 gives the following theorem.

**Theorem 4.4:** *Consider the system*

$$\begin{aligned} &\dot{x}_i = f_i(d, x_1, ..., x_i) + g_i(d, x_1, ..., x_i) x_{i+1} \quad i = 1, ..., n-1 \\ &\dot{x}_n = f_n(d,x) + g_n(d,x) u \\ &x = (x_1, ..., x_n)' \in \Re^n, d \in D, u \in \Re \end{aligned} \tag{4.13}$$

*where* $D \subset \Re^k$ *is a compact set ,* $f_i : D \times \Re^i \to \Re$, $g_i : D \times \Re^i \to \Re$ $(i = 1, ..., n)$ *are locally Lipschitz mappings with* $f_i(d,0,0) = 0$ *for all* $d \in D$ $(i = 1, ..., n)$, *which satisfy the following hypotheses:*

**(W5)** *There exist constants* $q, L, r > 0$ *such that* $|f_i(d,x)| \le \min\{q, L|x|\}$, $r \le g_i(d,x)$, *for all* $(d,x) \in D \times \Re^i$ $(i = 1, ..., n)$.

**(W6)** *There exists a constant* $R > 0$ *such that* $g_i(d,x) \le R$*, for all* $(d,x) \in D \times \Re^i$ *(*$i = 1,...,n-1$*).*

*Then there exist constants* $a_i, p_i \ge 0$ *(*$i = 1,...,n$*) such that the locally Lipschitz feedback law:*

$$u = \varphi_n(x) \tag{4.14}$$

*obtained by the recursive formula*

$$\varphi_{i+1}(x) := -a_{i+1} sat\big(p_{i+1}\big(x_{i+1} - \varphi_i(x)\big)\big),\ i = 1,...,n-1,\ \varphi_1(x_1) = -a_1 sat\big(p_1 x_1\big) \tag{4.15}$$

*robustly globally asymptotically stabilizes* $0 \in \Re^n$ *for system (4.13). Moreover, there exists a symmetric, positive definite matrix* $P \in \Re^{n \times n}$ *and a compact set* $\Omega \subset \Re^n$ *containing a neighborhood of* $0 \in \Re^n$ *such that the hypotheses (P1), (P2) of Theorem 2.2 hold with* $V(x) := x'Px$ *for appropriate* $T \in C^0(\Re^n;\Re^+)$*,* $G \in C^0(\Re^n;\Re^+)$*.*

**Sketch of proof:** By virtue of Lemma 4.3 it suffices to show that there exist $a_1, p_1 > 0$ such that hypotheses (W1-3) hold for the scalar subsystem:

$$\dot{x}_1 = f_1(d,x_1) + g(d,x_1)x_2$$

with $\varphi_1(x_1) = -a_1 sat\big(p_1 x_1\big)$, $P = [1]$, $S := \left\{x_1 \in \Re, |x_1| \le p_1^{-1}\right\}$ and appropriate $c > 0$. The proof of (W1) and (W3) is straightforward, while the proof of (W2) makes use of Lemma 2.5 (exactly as in the proof of Lemma 4.3). Details are left to the reader. $\lhd$

## 5. Concluding Remarks

The notion of relaxed Robust Control Lyapunov Function (RCLF) is introduced and is exploited for the design of robust feedback stabilizers for nonlinear systems. The development of the notion of the "relaxed" RCLF is important because even if a RCLF is known then the use of "relaxed" RCLF feedback design methodology usually results to different feedback designs from the ones obtained by the use of the standard RCLF design methodology; particularly, there is no need to make the derivative of RCLF negative everywhere. Moreover, in many cases "relaxed" RCLFs can be found, while RCLFs are not available. Consequently, the class of systems where Lyapunov-based feedback design principles can be applied is enlarged. Particularly, it is shown for systems with input constraints that "relaxed" RCLFs can be easily obtained, while RCLFs are not available. Moreover, it is shown that the use of "relaxed" RCLFs usually results to different feedback designs from the ones obtained by the use of the standard RCLF methodology. Using the "relaxed" RCLFs feedback design methodology, a simple controller that guarantees robust global stabilization of a perturbed chemostat model is provided. The proposed controller guarantees stabilization without assuming knowledge of the size of the uncertainty range.

## Appendix

**Proof of Lemma 2.4:** It suffices to show that $0 \in \Re^n$ is uniformly robustly Lyapunov stable. Let $s > 0$. Clearly, by virtue of the property of Uniform attractivity for bounded sets of initial states, there exists $T(s) \ge 0$ such that for every $(x_0, d) \in \Re^n \times M_D$ with $|x_0| \le s$ the solution $x(t; x_0, d)$ of (2.1) satisfies:

$$|x(t; x_0, d)| \le s \text{, for all } t \ge T(s) \tag{A1}$$

Let $r(s) := \max\{R(x_0) : |x_0| \le s\}$, where $R : \Re^n \to \Re^+$ is the continuous function involved in (2.2) and define $S := \{x \in \Re^n : |x| \le r(s)\}$ (a compact set) and $L(s) := \sup\left\{ \frac{(x-y)' P(F(d,x) - F(d,y))}{|x-y|^2} : x, y \in S, x \ne y \right\}$, where $P \in \Re^{n \times n}$ is the symmetric positive definite matrix involved in hypothesis (H3). Let $(x_0, d) \in \Re^n \times M_D$ with $|x_0| \le s$ and consider the solution $x(t; x_0, d)$ of (2.1). The evaluation of the derivative of the absolutely continuous function $V(t) = x'(t; x_0, d) P x(t; x_0, d)$, in conjunction with previous definitions, inequalities (A1), (2.2) and hypothesis (H3) gives:

$$\dot{V}(t) \le 2L(s)V(t) \text{, a.e. for } t \ge 0 \tag{A2}$$

Consequently, we obtain:

$$|x(t; x_0, d)| \le \sqrt{\frac{K_2}{K_1}} \exp(L(s)t) |x_0| \text{, } \forall t \ge t_0 \tag{A3}$$

where $K_1, K_2 > 0$ are constants satisfying $K_1 |x|^2 \le x'Px \le K_2 |x|^2$ for all $x \in \Re^n$. Combining (A1) and (A3) we conclude that for every $(x_0, d) \in \Re^n \times M_D$ with $|x_0| \le s$ the solution $x(t; x_0, d)$ of (2.1) satisfies the following estimate for all $\rho \ge s$:

$$|x(t; x_0, d)| \le \sqrt{\frac{K_2}{K_1}} \exp(L(\rho)T(\rho)) s \text{, } \forall t \ge 0 \tag{A4}$$

The above inequality guarantees Uniform Robust Lyapunov stability. Particularly, by virtue of (A4) it follows that for every $\varepsilon > 0$ there exists $\delta := \delta(\varepsilon) := \sqrt{\frac{K_1}{K_2}} \varepsilon \exp(-L(\varepsilon)T(\varepsilon)) > 0$ such that for all $(x_0, d) \in \Re^n \times M_D$ with $|x_0| \le \delta$ the solution $x(t; x_0, d)$ of (2.1) satisfies $|x(t)| \le \varepsilon$ for all $t \ge 0$. The proof is complete. ◁

**Proof of Theorem 2.2:** By virtue of Lemma 2.4 it suffices to show that there exists a continuous function $R : \Re^n \to \Re^+$ satisfying (2.2) and that the property of uniform attractivity for bounded sets of initial states holds. Standard arguments utilizing hypothesis (P2) (see [14,16]) and the fact that $x(t; x_0, d) \in \Omega$ for all $t \ge \hat{t}(x_0, d)$, show the existence of a function $\sigma \in KL$ such that for every $(x_0, d) \in \Re^n \times M_D$ the solution $x(t; x_0, d)$ of (2.1) satisfies $|x(t; x_0, d)| \le \sigma(|x(\hat{t}(x_0, d); x_0, d)|, t - \hat{t}(x_0, d))$ for all $t \ge \hat{t}(x_0, d)$. Using the previous estimate and hypothesis (P1) we obtain:

$$|x(t; x_0, d)| \le \sigma(G(x_0), t - \hat{t}(x_0, d)) \text{, for all } t \ge \hat{t}(x_0, d) \tag{A5}$$

$$|x(t; x_0, d)| \le \max\{G(x_0), \sigma(G(x_0), 0)\} \text{, for all } t \ge 0 \tag{A6}$$

Inequality (A6) shows that the continuous function $R(x_0):=\max\left\{G(x_0)\,,\sigma\left(G(x_0)\,,0\right)\right\}$ satisfies inequality (2.2). Moreover, inequality (A5) and the fact $\hat{t}(x_0,d)\in[0,T(x_0)]$ show that for every $\varepsilon>0$, $s\geq 0$, $(x_0,d)\in\Re^n\times M_D$ with $|x_0|\leq s$ the solution $x(t;x_0,d)$ of (2.1) satisfies $|x(t;x_0,d)|\leq\varepsilon$ for all $t\geq T(\varepsilon,s)$, where $T(\varepsilon,s):=g(\varepsilon,s)+\tau(s)$, $\tau(s):=\max\left\{T(x_0):|x_0|\leq s\right\}$ and $g(\varepsilon,s)\geq 0$ is any time satisfying $\sigma\left(r(s),g(\varepsilon,s)\right)\leq\varepsilon$ with $r(s):=\max\left\{G(x_0):|x_0|\leq s\right\}$. The proof is complete. ◁

**Proof of Lemma 2.5:** First notice that inequalities (2.3a,b), (2.4) in conjunction with Corollary 8.2 in [5] imply that the following implications hold for every $(x_0,d)\in\Re^n\times M_D$ :

If $h(x(t;x_0,d))\in(0,b)$, $\dot{x}(t;x_0,d)=F(d(t),x(t;x_0,d))$ and $\frac{d}{dt}h(x(t;x_0,d))$ exists then $\frac{d}{dt}h(x(t;x_0,d))\leq 0$ (A7)

If $h(x(t;x_0,d))>0$, $\dot{x}(t;x_0,d)=F(d(t),x(t;x_0,d))$ and $\frac{d}{dt}q(x(t;x_0,d))$ exists then

$$\frac{d}{dt}q(x(t;x_0,d))\leq-\delta\left(h(x(t;x_0,d))\right) \tag{A8}$$

If $h(x(t;x_0,d))>0$, $\dot{x}(t;x_0,d)=F(d(t),x(t;x_0,d))$ and $\frac{d}{dt}W(x(t;x_0,d))$ exists then

$$\frac{d}{dt}W(x(t;x_0,d))\leq KW(x(t;x_0,d)) \tag{A9}$$

where $q(x(t;x_0,d)):=h(x(t;x_0,d))-a(x(t;x_0,d))$.

Let $\hat{\varepsilon}\in(0,b)$. Notice that implication (A7) guarantees that the set $\Omega:=\left\{x\in\Re^n : h(x)\leq\hat{\varepsilon}\right\}$ is positively invariant for system (2.1). Indeed, if $x_0\in\Omega$ then for every $d\in M_D$ it holds that $x(t;x_0,d)\in\Omega$ for all $t\in[0,t_{\max})$, where $t_{\max}=t_{\max}(x_0,d)$ is the maximal existence time of the solution. In order to show positive invariance of $\Omega:=\left\{x\in\Re^n : h(x)\leq\hat{\varepsilon}\right\}$, we use the following contradiction argument: suppose that there exists $x_0\in\Omega$, $d\in M_D$ and $t\in(0,t_{\max})$ such that $x(t;x_0,d)\notin\Omega$, i.e., $h(x(t;x_0,d))>\hat{\varepsilon}$. Exploiting continuity of the mapping $\tau\to h(x(\tau;x_0,d))$ we guarantee that the set $A:=\{\tau\in[0,t]:h(x(\tau;x_0,d))=\hat{\varepsilon}\}$ is non-empty. Let $T:=\sup A$. Notice that continuity of the mapping $\tau\to h(x(\tau;x_0,d))$ implies $T<t$, $h(x(T;x_0,d))=\hat{\varepsilon}$ and $h(x(\tau;x_0,d))>\hat{\varepsilon}$ for all $\tau\in(T,t]$. Without loss of generality we may also assume that $h(x(\tau;x_0,d))<b$ for all $\tau\in(T,t)$ (possibly by replacing $t$ by $\hat{t}:=\inf B$, where $B:=\{\tau\in[T,t]:h(x(\tau;x_0,d))\geq b\}$). Absolute continuity of the mapping $\tau\to h(x(\tau;x_0,d))$ implies that $\hat{\varepsilon}<h(x(t;x_0,d))=h(x(T;x_0,d))+\int_T^t\frac{d}{d\tau}h(x(\tau;x_0,d))\,d\tau\leq\hat{\varepsilon}$, a contradiction.

Next, we consider the case $x_0\notin\Omega$. Let arbitrary $x_0\in\Re^n$, with $h(x_0)>\hat{\varepsilon}$, $d\in M_D$ and consider the solution $x(t;x_0,d)$ of (2.1). Define the set $\{t\geq 0:x(t;x_0,d)\notin\Omega\}$. Clearly this set is non-empty (since $0\in\{t\geq 0:x(t;x_0,d)\notin\Omega\}$). We next claim that $\hat{t}(x_0,d):=\sup\{t\geq 0:x(t;x_0,d)\notin\Omega\}\leq\frac{q(x_0)+L-\hat{\varepsilon}}{\min\left\{\delta(s):\hat{\varepsilon}\leq s\leq L+q(x_0)\right\}}$, where $q(x_0):=h(x_0)-a(x_0)$, $L:=\sup_{x\in\Re^n}a(x)$. The proof is made by contradiction. Suppose that this is not the case. Then there exists $t>\frac{q(x_0)+L-\hat{\varepsilon}}{\min\left\{\delta(s):\hat{\varepsilon}\leq s\leq L+q(x_0)\right\}}$ with $h(x(t;x_0,d))>\hat{\varepsilon}$. Since $\Omega:=\left\{x\in\Re^n : h(x)\leq\hat{\varepsilon}\right\}$ is positively invariant for system (2.1), this implies that $h(x(\tau;x_0,d))>\hat{\varepsilon}$ for all $\tau\in[0,t]$. Consequently, it follows from (A8) that $\frac{d}{d\tau}q(x(\tau;x_0,d))\leq 0$, a.e. on $[0,t]$, where $q(x(t;x_0,d)):=h(x(t;x_0,d))-a(x(t;x_0,d))$. Therefore the mapping $[0,t]\ni\tau\to q(x(\tau;x_0,d))$ is non-increasing, i.e., $q(x(\tau;x_0,d))\leq q(x_0)$ for all $\tau\in[0,t]$. Thus, it holds that $\hat{\varepsilon} \le h(x(\tau;x_0,d)) \le q(x_0)+L$ for all $\tau \in [0,t]$, where $L := \sup_{x\in\Re^n} a(x)$. Differential inequality (A8) and the fact that $\hat{\varepsilon} \le h(x(\tau;x_0,d)) \le q(x_0)+L$ for all $\tau \in [0,t]$ gives $\frac{d}{d\tau} q(x(\tau;x_0,d)) \le -\min\{\delta(s): \hat{\varepsilon} \le s \le L+q(x_0)\}$ a.e. on $[0,t]$. Thus we obtain $q(x(t;x_0,d)) \le q(x_0) - t\min\{\delta(s): \hat{\varepsilon} \le s \le L+q(x_0)\}$, which directly implies $h(x(t;x_0,d)) \le q(x_0)+L-t\min\{\delta(s): \hat{\varepsilon} \le s \le L+q(x_0)\}$. The latter inequality combined with the hypothesis $t > \frac{q(x_0)+L-\hat{\varepsilon}}{\min\{\delta(s): \hat{\varepsilon} \le s \le L+q(x_0)\}}$ gives $h(x(t;x_0,d)) \le \hat{\varepsilon}$, a contradiction.

Since $\hat{t}(x_0,d) := \sup\{t \ge 0 : x(t;x_0,d) \notin \Omega\} \le \frac{q(x_0)+L-\hat{\varepsilon}}{\min\{\delta(s): \hat{\varepsilon} \le s \le L+q(x_0)\}}$ and $\hat{t}(x_0,d) > 0$, it follows from (A9) and previous definitions that $\frac{d}{d\tau} W(x(\tau;x_0,d)) \le KW(x(\tau;x_0,d))$, a.e. on $[0,\hat{t}(x_0,d))$. The previous differential inequality implies $W(x(\tau;x_0,d)) \le \exp\left(K\frac{q(x_0)+L-\hat{\varepsilon}}{\min\{\delta(s): \hat{\varepsilon} \le s \le L+q(x_0)\}}\right)W(x_0)$ for all $\tau \in [0,\hat{t}(x_0,d))$. Since $W$ is radially unbounded, it follows that the solution of (2.1) exists on $[0,\hat{t}(x_0,d)]$ and satisfies $x(\hat{t}(x_0,d);x_0,d) \in \Omega$, $W(x(\tau;x_0,d)) \le \exp\left(K\frac{q(x_0)+L-\hat{\varepsilon}}{\min\{\delta(s): \hat{\varepsilon} \le s \le L+q(x_0)\}}\right)W(x_0)$ for all $\tau \in [0,\hat{t}(x_0,d)]$. Consequently, $t_{\max} > \hat{t}(x_0,d)$ and since the set $\Omega := \{x \in \Re^n : h(x) \le \hat{\varepsilon}\}$ is positively invariant for system (2.1) it follows that $x(t;x_0,d) \in \Omega$ for all $t \in [\hat{t}(x_0,d), t_{\max})$

Since $W$ is radially unbounded, it follows that for all $s \ge \min_{x\in\Re^n} W(x)$ the set $\{x \in \Re^n : W(x) \le s\}$ is non-empty and compact. Define $r(s) := \max\{|x| : x \in \Re^n, W(x) \le s\}$ for $s \ge \min_{x\in\Re^n} W(x)$, which is a non-decreasing function. Define the continuous function $B(x_0) := \exp\left(K\frac{\max\{0, q(x_0)+L-\hat{\varepsilon}\}}{\min\{\delta(s): \hat{\varepsilon} \le s \le \hat{\varepsilon}+\max\{0, q(x_0)+L\}\}}\right)W(x_0)$ and let $R(x_0) \ge r(B(x_0))$. By distinguishing the cases $x_0 \in \Omega$ (which implies $\hat{t}(x_0,d)=0$) and $x_0 \notin \Omega$, we notice that property (P1) of Theorem 2.2 holds with $\Omega := \{x \in \Re^n : h(x) \le \hat{\varepsilon}\}$ and functions $T \in C^0(\Re^n;\Re^+)$, $G \in C^0(\Re^n;\Re^+)$ defined as follows:

$$T(x_0) := \frac{\max\{0, q(x_0)+L-\hat{\varepsilon}\}}{\min\{\delta(s): \hat{\varepsilon} \le s \le \hat{\varepsilon}+\max\{0, q(x_0)+L\}\}},\quad G(x_0) := \tilde{r}(B(x_0))$$

where $\tilde{r}(s) := \int_s^{s+1} r(w)dw$. The proof is complete. $\lhd$

**Proof of Theorem 2.6:** Without loss of generality and since $h(0) < 0$ we may assume that the neighborhood $N$ involved in hypothesis (R4) satisfies $N \subset \{x \in \Re^n : h(x) < 0\}$. Let $r > 0$ with $\{x \in \Re^n : |x| \le 2r\} \subset N$.

The construction of the feedback will be accomplished in three steps:

Step 1: Construction of preliminary feedback laws, which work on certain sets of the state space.

Step 2: Definition of the required feedback law by patching together the preliminary feedback laws of the previous step.

Step 3: Proof of URGAS by using Theorem 2.2 and Lemma 2.5.

Step 1: Construction of preliminary feedback laws, which work on certain sets of the state space.

Using hypothesis (R1), convexity of $U \subseteq \Re^m$ and standard partition of unity arguments, we construct a smooth feedback law $k_1 : \{x \in \Re^n : h(x) > 0\} \to U$ such that the following inequalities hold for all $x \in \{x \in \Re^n : h(x) > 0\}$:

$$\sup_{d \in D} \nabla h(x)\big(f(d,x) + g(d,x)k_1(x)\big) \le -\frac{1}{2}\delta(h(x)) \tag{A10}$$

$$\sup_{d \in D} \nabla W(x)\big(f(d,x) + g(d,x)k_1(x)\big) \le KW(x) + 1 \tag{A11}$$

Using hypothesis (R2), convexity of $U \subseteq \Re^m$ and standard partition of unity arguments, we construct a smooth feedback law $k_2 : \{x \in \Re^n : h(x) < \varepsilon, x \neq 0\} \to U$ such that the following inequality holds for all $x \in \{x \in \Re^n : h(x) < \varepsilon, x \neq 0\}$:

$$\sup_{d \in D} \nabla V(x)\big(f(d,x) + g(d,x)k_2(x)\big) < 0 \tag{A12}$$

Using hypothesis (R3), convexity of $U \subseteq \Re^m$ and standard partition of unity arguments, we construct a smooth feedback law $k_3 : \{x \in \Re^n : 0 < h(x) < \varepsilon\} \to U$ such that the following inequalities hold for all $x \in \{x \in \Re^n : 0 < h(x) < \varepsilon\}$

$$\sup_{d \in D} \nabla h(x)\big(f(d,x) + g(d,x)k_3(x)\big) \le -\frac{1}{2}\delta(h(x)) \tag{A13}$$

$$\sup_{d \in D} \nabla W(x)\big(f(d,x) + g(d,x)k_3(x)\big) \le KW(x) + 1 \tag{A14}$$

$$\sup_{d \in D} \nabla V(x)\big(f(d,x) + g(d,x)k_3(x)\big) < 0 \tag{A15}$$

<u>Step 2:</u> Definition of the required feedback law by patching together the preliminary feedback laws of the previous step.

Let $p : \Re \to [0,1]$ a smooth non-decreasing function with $p(x) = 0$ for all $x \le 0$ and $p(x) = 1$ for all $x \ge 1$. We define:

$$k(x) := k_1(x) \text{ for all } x \in \left\{x \in \Re^n : h(x) > \frac{4\varepsilon}{5}\right\} \tag{A16}$$

$$k(x) := k_2(x) \text{ for all } x \in \left\{x \in \Re^n : h(x) < \frac{\varepsilon}{5}\right\} \cap \left\{x \in \Re^n : |x| > 2r\right\} \tag{A17}$$

$$k(x) := k_3(x) \text{ for all } x \in \left\{x \in \Re^n : \frac{2\varepsilon}{5} < h(x) < \frac{3\varepsilon}{5}\right\} \tag{A18}$$

$$k(x) := \left(1 - p\left(\frac{5h(x)}{\varepsilon} - 3\right)\right)k_3(x) + p\left(\frac{5h(x)}{\varepsilon} - 3\right)k_1(x) \text{ for all } x \in \left\{x \in \Re^n : \frac{3\varepsilon}{5} \le h(x) \le \frac{4\varepsilon}{5}\right\} \tag{A19}$$

$$k(x) := \left(1 - p\left(\frac{5h(x)}{\varepsilon} - 1\right)\right)k_2(x) + p\left(\frac{5h(x)}{\varepsilon} - 1\right)k_3(x) \text{ for all } x \in \left\{x \in \Re^n : \frac{\varepsilon}{5} \le h(x) \le \frac{2\varepsilon}{5}\right\} \tag{A20}$$

$$k(x) := \tilde{k}(x) \text{ for all } x \in \left\{x \in \Re^n : |x| < r\right\} \tag{A21}$$

$$k(x) := \left(1 - p\left(\frac{|x|^2 - r^2}{3r^2}\right)\right)\tilde{k}(x) + p\left(\frac{|x|^2 - r^2}{3r^2}\right)k_2(x) \text{ for all } x \in \left\{x \in \Re^n : r \le |x| \le 2r\right\} \tag{A22}$$

where $\tilde{k} : N \to U$ is the locally Lipschitz mapping involved in hypothesis (R4). Convexity of $U \subseteq \Re^m$ and the above definitions imply that $k(x) \in U$ for all $x \in \Re^n$. Notice that the mapping $k : \Re^n \to U$ defined above is locally Lipschitz with $k(0) = 0$. By virtue of (A10), (A11), (A12), (A13), (A14), (A15) and definitions (A16-22) we obtain the following inequalities:

$$\sup_{d \in D} \nabla h(x)\big(f(d,x) + g(d,x)k(x)\big) \le -\frac{1}{2}\delta(h(x)) \text{, for all } x \in \left\{x \in \Re^n : h(x) \ge \frac{2\varepsilon}{5}\right\} \tag{A23}$$

$$\sup_{d \in D} \nabla W(x)\big(f(d,x) + g(d,x)k(x)\big) \le KW(x) + 1 \text{, for all } x \in \left\{x \in \Re^n : h(x) \ge \frac{2\varepsilon}{5}\right\} \tag{A24}$$

$$\sup_{d \in D} \nabla V(x)\big(f(d,x) + g(d,x)k(x)\big) < 0 \text{, for all } x \in \left\{x \in \Re^n : x \ne 0, h(x) \le \frac{3\varepsilon}{5}\right\} \tag{A25}$$

<u>Step 3:</u> Proof of URGAS by using Theorem 2.2 and Lemma 2.5.

First notice that by virtue of hypotheses (Q1-3) and the fact that the mapping $k : \Re^n \to U$ defined above is locally Lipschitz with $k(0) = 0$, it follows that the closed-loop system (1.1) with $u = k(x)$ satisfies hypotheses (H1-3). Next define $\tilde{h}(x) := h(x) - \frac{2\varepsilon}{5}$, $\tilde{\delta}(s) := \frac{1}{2}\delta\left(s + \frac{2\varepsilon}{5}\right)$, $\tilde{K} := K + 1$, $\tilde{W}(x) := W(x) + 1$ and. It should be noticed that by virtue of inequalities (A23), (A24), all requirements of Lemma 2.5 hold with $F(d,x) := f(d,x) + g(d,x)k(x)$, $a(x) \equiv 0$, arbitrary $b > \frac{\varepsilon}{5}$, $\hat{\varepsilon} := \frac{\varepsilon}{5}$ and $\tilde{h}, \tilde{\delta}, \tilde{K}, \tilde{W}$ in place of $h, \delta, K, W$, respectively. Therefore, there exist functions $T \in C^0(\Re^n;\Re^+)$, $G \in C^0(\Re^n;\Re^+)$ such that property (P1) of Theorem 2.2 holds with $\Omega := \left\{x \in \Re^n : \tilde{h}(x) \le \tilde{\varepsilon}\right\} = \left\{x \in \Re^n : h(x) \le \frac{3\varepsilon}{5}\right\}$. On the other hand, inequality (A25) guarantees that property (P2) of Theorem 2.2 holds with $\Omega := \left\{x \in \Re^n : h(x) \le \frac{3\varepsilon}{5}\right\}$. Consequently, we may conclude from Theorem 2.2 that $0 \in \Re^n$ is URGAS for the closed-loop system (1.1) with $u = k(x)$. The proof is complete. $\triangleleft$

**Proof of Lemma 4.3:** Let $p > 0$ sufficiently large so that:

$$p \ge c^{-1} \tag{A26}$$

$$p\big(\tilde{c}r + q + b\big) \ge \frac{\mu}{2} + L + C|k| + \frac{1}{2\mu K}\big(C|P| + C|k| + L + L|k|\big)^2 \tag{A27}$$

where $K := \min_{|x|=1} x'Px$.

We start by proving the analogue of (W1) for system (4.9), (4.12). First notice that the following equality holds for all $(d,x,y) \in D \times \Re^n \times \Re$

$$\tilde{x}'\tilde{P}\tilde{F}(d,\tilde{x},\tilde{k}\,\tilde{x}) = x'PF(d,x,y) + (y - k'x)(f(d,x,y) - \tilde{a}pg(d,x,y)(y - k'x) - k'F(d,x,y))$$

Exploiting (4.10) we get from the above equation for all $(d,x,y)\in D\times\Re^{n}\times\Re$ :

$$\tilde{x}'\tilde{P}\tilde{F}(d,\tilde{x},\tilde{k}\,\tilde{x})\le -\mu x'Px + x'P(F(d,x,y)-F(d,x,k'x)) + (y-k'x)(f(d,x,y)-\tilde{a}pg(d,x,y)(y-k'x)-k'F(d,x,y))$$

Hypotheses (W3), (W4) in conjunction with the above inequality imply for all $(d,x,y)\in D\times\Re^{n}\times\Re$ :

$$\tilde{x}'\tilde{P}\tilde{F}(d,\tilde{x},\tilde{k}\,\tilde{x})\le -\mu x'Px + (C|P|+C|k|+L+L|k|)\,|x|\,|y-k'x| - (\tilde{a}pr-L-C|k|)|y-k'x|^{2}$$

Completing the squares we get for all $(d,x,y)\in D\times S\times\Re$ :

$$\tilde{x}'\tilde{P}\tilde{F}(d,\tilde{x},\tilde{k}\,\tilde{x})\le -\frac{1}{2}\mu x'Px + \frac{1}{2\mu K}(C\gamma|P|+C|k|+L+L|k|)^{2}\ |y-k'x|^{2} - (\tilde{a}pr-L-C\gamma|k|)|y-k'x|^{2}$$

where $K:=\min_{|x|=1} x'Px$ . Finally, notice that the above inequality in conjunction with the facts $\tilde{x}'\tilde{P}\tilde{x}=x'Px+(y-k'x)^{2}$ , $\tilde{a}>\tilde{c}+r^{-1}q$ and (A27) implies for all $(d,x,y)\in D\times S\times\Re$ :

$$\tilde{x}'\tilde{P}\tilde{F}(d,\tilde{x},\tilde{k}\,\tilde{x})\le -\frac{\mu}{2}\tilde{x}'\tilde{P}\tilde{x} \tag{A28}$$

Moreover, the equality $\tilde{\varphi}(\tilde{x}):=-\tilde{a}\,sat(p(y-\varphi(x)))=\tilde{k}\,\tilde{x}$ for all $\tilde{x}\in\tilde{S}:=\{(x,y)\in S\times\Re:|y-\varphi(x)|\le p^{-1}\}\subset\Re^{n+1}$ holds automatically, by virtue of (4.11). Therefore, the analogue of hypothesis (W1) holds for system (4.9), (4.12).

We continue by showing the analogue of hypothesis (W2). Let $(x_0,y_0,d,v)\in\Re^{n}\times\Re\times M_D\times M_{[-\tilde{c},\tilde{c}]}$ the solution $x(t)$ of (4.9), (4.12) with $u=-\tilde{a}\,sat(p(y-\varphi(x)))+v$ , $(x(0),y(0))=(x_0,y_0)$ corresponding to inputs $(d,v)\in M_D\times M_{[-\tilde{c},\tilde{c}]}$ . Hypothesis (W3) implies that

$$|\dot{x}(t)|\le C|x(t)|+C+C|y(t)-\varphi(x(t))|\text{, a.e. for } t\in[0,t_{\max}) \tag{A29}$$

where $t_{\max}>0$ is the maximal existence time of the solution. Define the absolutely function $Y(t)=\max\{a+p^{-1},|y(t)|\}$. Notice that if $y(t)>a+p^{-1}$ then $y(t)-\varphi(x(t))>p^{-1}$ and consequently $u(t)=-\tilde{a}+v(t)$ . Moreover, if $\dot{y}(t)=f(d,x(t),y(t))+g(d,x(t),y(t))u(t)$ , then the inequalities $\tilde{a}>\tilde{c}+r^{-1}q$ and $v(t)\le\tilde{c}$ in conjunction with hypothesis (W4) imply that $\dot{y}(t)\le 0$ . Similarly, the case $y(t)<-a-p^{-1}$ implies $\dot{y}(t)\ge 0$ . Therefore, we obtain:

$$\dot{Y}(t)\le 0\text{, a.e. for } t\in[0,t_{\max}) \tag{A30}$$

The above differential inequality implies that the mapping $t\to Y(t)$ is non-increasing and consequently we obtain the estimate:

$$|y(t)|\le a+p^{-1}+|y_0|\text{, for all } t\in[0,t_{\max}) \tag{A31}$$

Since $|\varphi(x(t))|\le a$ , we get from (A29) and (A31):

$$|\dot{x}(t)|\le C|x(t)|+C(1+2a+p^{-1}+|y_0|)\text{, a.e. for } t\in[0,t_{\max}) \tag{A32}$$

which directly implies

$$|x(t)|\le(|x_0|+1+2a+p^{-1}+|y_0|)\exp(Ct)\text{ , for all } t\in[0,t_{\max}) \tag{A33}$$

Inequalities (A31) and (A33) guarantee that $t_{\max} = +\infty$.

We are now in a position to apply Lemma 2.5 with $h(x,y) := y - \varphi(x) - \lambda p^{-1}$, $a(x,y) = -\varphi(x)$, $W(x,y) = |x| + \max\{a + p^{-1}, |y(t)|\} + 1$, $b := (1-\lambda)p^{-1}$, $\delta(h(x,y)) := r(\lambda\tilde{a} - \tilde{c}) - q$, where $\frac{\tilde{c}r + q + b}{\tilde{a}r} < \lambda < 1$. Using (A29), (A30), hypotheses (W3), (W4), inequalities (A26), $|\varphi(x)| \le a$ and definition $\tilde{\varphi}(\tilde{x}) := -\tilde{a}\, sat(p(y - \varphi(x)))$ with $\tilde{a} > \tilde{c} + r^{-1}(q+b)$, it may be shown that:

$$\sup_{d \in D, v \in [-\tilde{c},\tilde{c}]} \nabla h(x,y) \tilde{F}(d,x,y,\tilde{\varphi}(x,y)+v) \le 0 \text{, for almost all } x \in \Re^n \text{ with } 0 < h(x) < \varepsilon \quad \text{(A34)}$$

$$\sup_{d \in D, v \in [-\tilde{c},\tilde{c}]} (\nabla h(x,y) - \nabla a(x,y)) \tilde{F}(d,x,y,\tilde{\varphi}(x,y)+v) \le -\delta(h(x,y)) \text{, for almost all } x \in \Re^n \text{ with } h(x) > 0 \quad \text{(A35)}$$

$$\sup_{d \in D, v \in [-\tilde{c},\tilde{c}]} \nabla W(x,y) F(d,x,y,\tilde{\varphi}(x,y)+v) \le KW(x,y) \text{, for almost all } x \in \Re^n \text{ with } h(x) > 0 \quad \text{(A36)}$$

with $K := C(1+a)$. Therefore, there exist mappings $T_1 \in C^0(\Re^n \times \Re; \Re^+)$, $Q_1 \in C^0(\Re^n \times \Re; \Re^+)$ such that for every $(x_0, y_0, d, v) \in \Re^n \times \Re \times M_D \times M_{[-\tilde{c},\tilde{c}]}$, there exists $\hat{t}_1(x_0,y_0,d,v) \in [0, T_1(x_0,y_0)]$ in such a way that the solution $(x(t), y(t))$ of (4.9), (4.12) with $u = -\tilde{a}\, sat(p(y - \varphi(x))) + v$, $(x(0), y(0)) = (x_0, y_0)$ corresponding to inputs $(d,v) \in M_D \times M_{[-\tilde{c},\tilde{c}]}$ satisfies $y(t) - \varphi(x(t)) \le p^{-1}$ for all $t \ge \hat{t}_1(x_0,y_0,d,v)$ and $|(x(t),y(t))| \le Q_1(x_0,y_0)$ for all $t \in [0, \hat{t}_1(x_0,y_0,d,v)]$.

Applying Lemma 2.5 (again) with $h(x,y) := \varphi(x) - y - \lambda p^{-1}$, $a(x,y) = \varphi(x)$, $W(x,y) = |x| + \max\{a + p^{-1}, |y(t)|\} + 1$, $b := (1-\lambda)p^{-1}$, $\delta(h(x,y)) := r(\lambda\tilde{a} - \tilde{c}) - q$, where $\frac{\tilde{c}r + q + b}{\tilde{a}r} < \lambda < 1$. Using (again) (A29), (A30), hypotheses (W3), (W4), inequalities (A26), $|\varphi(x)| \le a$ and definition $\tilde{\varphi}(\tilde{x}) := -\tilde{a}\, sat(p(y - \varphi(x)))$ with $\tilde{a} > \tilde{c} + r^{-1}(q+b)$, it may be shown that inequalities (A34), (A35), (A36) hold as well with $K := C(1+a)$. Therefore, there exist mappings $T_2 \in C^0(\Re^n \times \Re; \Re^+)$, $Q_2 \in C^0(\Re^n \times \Re; \Re^+)$ such that for every $(x_0, y_0, d, v) \in \Re^n \times \Re \times M_D \times M_{[-\tilde{c},\tilde{c}]}$, there exists $\hat{t}_2(x_0,y_0,d,v) \in [0, T_2(x_0,y_0)]$ in such a way that the solution $(x(t), y(t))$ of (4.9), (4.12) with $u = -\tilde{a}\, sat(p(y - \varphi(x))) + v$, $(x(0), y(0)) = (x_0, y_0)$ corresponding to inputs $(d,v) \in M_D \times M_{[-\tilde{c},\tilde{c}]}$ satisfies $y(t) - \varphi(x(t)) \ge -p^{-1}$ for all $t \ge \hat{t}_2(x_0,y_0,d,v)$ and $|(x(t),y(t))| \le Q_2(x_0,y_0)$ for all $t \in [0, \hat{t}_2(x_0,y_0,d,v)]$.

We conclude that for every $(x_0, y_0, d, v) \in \Re^n \times \Re \times M_D \times M_{[-\tilde{c},\tilde{c}]}$, there exists $\bar{t}(x_0,y_0,d,v) \in [0, \hat{T}(x_0,y_0)]$, where $\hat{T}(x_0,y_0) := \max\{T_1(x_0,y_0), T_2(x_0,y_0)\}$, $\bar{t}(x_0,y_0,d,v) := \max\{\hat{t}_1(x_0,y_0,d,v), \hat{t}_2(x_0,y_0,d,v)\}$ in such a way that the solution $(x(t), y(t))$ of (4.9), (4.12) with $u = -\tilde{a}\, sat(p(y - \varphi(x))) + v$, $(x(0), y(0)) = (x_0, y_0)$ corresponding to inputs $(d,v) \in M_D \times M_{[-\tilde{c},\tilde{c}]}$ satisfies $|y(t) - \varphi(x(t))| \le p^{-1}$ for all $t \ge \bar{t}(x_0,y_0,d,v)$ and $|(x(t),y(t))| \le \hat{Q}(x_0,y_0)$ for all $t \in [0, \bar{t}(x_0,y_0,d,v)]$, where $\hat{Q}(x_0,y_0) := \max\{\hat{Q}_1(x_0,y_0), \hat{Q}_2(x_0,y_0)\}$.

Since $p^{-1} \le c$ (a consequence of (A26)), hypothesis (W2) and inequalities (A31), (A33) imply that for every $(x_0, y_0, d, v) \in \Re^n \times \Re \times M_D \times M_{[-\tilde{c},\tilde{c}]}$, there exists $\tilde{t}(x_0,y_0,d,v) \in [0, \tilde{T}(x_0,y_0)]$, where

$$\tilde{T}(x_0,y_0) := \hat{T}(x_0,y_0) + \max\left\{ T(x) : |x| \le \left(|x_0| + 1 + 2a + p^{-1} + |y_0|\right)\exp(C\hat{T}(x_0,y_0)) \right\}$$

in such a way that the solution $(x(t), y(t))$ of (4.9), (4.12) with $u = -\tilde{a}\, sat(p(y - \varphi(x))) + v$, $(x(0), y(0)) = (x_0, y_0)$ corresponding to inputs $(d,v) \in M_D \times M_{[-\tilde{c},\tilde{c}]}$ satisfies $(x(t), y(t)) \in \tilde{S}$ for all $t \ge \tilde{t}(x_0,y_0,d,v)$ and $|(x(t),y(t))| \le \tilde{Q}(x_0,y_0)$ for all $t \in [0, \tilde{t}(x_0,y_0,d,v)]$, where

$$\widetilde{Q}(x_0,y_0):=\hat{Q}(x_0,y_0)+\max\left\{Q(x):|x|\le\left(|x_0|+1+2a+p^{-1}+|y_0|\right)\exp(C\hat{T}(x_0,y_0))\right\}+a+p^{-1}+|y_0|$$

Finally, we have:

$$\nabla\widetilde{\varphi}(\widetilde{x})=0 \text{ , provided that } |y-\varphi(x)|>p^{-1}$$

$$\nabla\widetilde{\varphi}(\widetilde{x})=-\widetilde{a}p\left[-\nabla\varphi(x),1\right] \text{, a.e. for } (x,y)\in\left\{(x,y)\in\Re^{n+1}:|y-\varphi(x)|<p^{-1}\right\}$$

Consequently, if there exists constant $R>0$ such that $g(d,x,y)\le R$, for all $(d,x,y)\in D\times\Re^n\times\Re$ then we obtain $\nabla\widetilde{\varphi}(\widetilde{x})\widetilde{F}(d,\widetilde{x},\widetilde{\varphi}(\widetilde{x})+v)=\widetilde{a}p\nabla\varphi(x)F(d,x,y)-\widetilde{a}p\left(f(d,x,y)+g(d,x,y)\gamma(\widetilde{\varphi}(\widetilde{x})+v)\right)$, a.e. for $(x,y)\in\left\{(x,y)\in\Re^{n+1}:|y-\varphi(x)|<p^{-1}\right\}$ and the inequality $\left|\nabla\widetilde{\varphi}(\widetilde{x})\widetilde{F}(d,\widetilde{x},\widetilde{\varphi}(\widetilde{x})+v)\right|\le\widetilde{b}$, for almost all $\widetilde{x}\in\Re^n$ and all $(d,v)\in D\times[-\widetilde{c},\widetilde{c}]$, where $\widetilde{b}:=\widetilde{a}p\left(q+R\widetilde{a}+R\widetilde{c}+b\right)$, follows directly from hypotheses (W3), (W4). The proof is complete. ◁